\documentstyle[twoside,12pt]{article}

\newtheorem{prop}{Proposition}[subsection]
\newtheorem{opr}[prop]{Definition}
\newtheorem{theo}[prop]{Theorem}

\newtheorem{rem}[prop]{Remark}
\newtheorem{coro}[prop]{Corollary}
\newtheorem{exam}[prop]{Example}

\newtheorem{lem}[prop]{Lemma}

\newcommand{\proof}{{\em Proof. }}

\title{Dual Polyhedra and  Mirror Symmetry for
Calabi-Yau Hypersurfaces in Toric Varieties}

\author{Victor V. Batyrev \\
 Universit\"at-GH-Essen, Fachbereich  6,  Mathematik \\
 Universit\"atsstr. 3,  45141  Essen  \\
Federal Republic of Germany \\
e-mail: matf$\emptyset \emptyset$@vm.hrz.uni-essen.de}

\begin{document}

\date{}

\maketitle

\footnote{Supported by DFG, Forschungsschwerpunkt  Komplexe
Mannigfaltigkeiten.}

\begin{abstract}
We  consider   families ${\cal F}(\Delta)$ consisting
of complex $(n-1)$-dimensional projective algebraic
compactifications of $\Delta$-regular affine
hypersurfaces $Z_f$ defined by Laurent polynomials $f$ with a
fixed $n$-dimensional Newton polyhedron $\Delta$ in
$n$-dimensional algebraic torus ${\bf T} = ({\bf C}^*)^n$.
If the  family  ${\cal F}(\Delta)$ defined by a
Newton polyhedron  $\Delta$ consists  of
$(n-1)$-dimensional Calabi-Yau varieties then
the dual, or polar,  polyhedron $\Delta^*$ in the dual space
defines  another
family ${\cal F}(\Delta^*)$ of Calabi-Yau varieties, so that we
obtain the remarkable duality between two {\em different families} of
Calabi-Yau varieties.
It is shown  that the properties of this duality
coincide with the properties of
{\em Mirror Symmetry}  discovered by physicists
for Calabi-Yau $3$-folds. Our method allows to construct
many new examples of Calabi-Yau $3$-folds and new candidates for their
mirrors which  were previously unknown for physicists.

We conjecture that there exists an isomorphism between
two conformal field theories corresponding  to Calabi-Yau varieties
from two families ${\cal F}(\Delta)$ and ${\cal F}(\Delta^*)$.
\end{abstract}

\section{Introduction}
\bigskip

\hspace*{\parindent}
Calabi-Yau  $3$-folds   caught much attention from
theoretical physics because of their connection with the superstring
theory \cite{witt}. Physicists discovered a duality  for
Calabi-Yau $3$-folds which is  called {\em Mirror Symmetry}
\cite{aspin1,aspin2,aspin3,cand.sch,dixon,gree,lerche,lynker}. This duality
defines a  correspondence between two topologically different
Calabi-Yau $3$-folds
$V$ and $V'$ such that the Hodge numbers of $V$ and $V'$ satisfy the
relations
\begin{equation}
  h^{1,1}(V) = h^{2,1}(V'), \;\; h^{1,1}(V') = h^{2,1}(V).
\label{main.relation}
\end{equation}
Thousands examples of Calabi-Yau $3$-folds obtained from
hypersurfaces in weighted projective spaces have shown  a
a striking symmetry for possible  pairs of integers
$(h^{1,1}, h^{2,1})$ relative to the transposition interchanging
$h^{1,1}$ and $h^{2,1}$ \cite{cand.sch,kreuzer,schimm}.
This fact gives an empirical evidence
in favor of the conjectural duality. On the other hand,
physicists have proposed some explicit constructions
of mirror manifolds \cite{berglund,gree} for several classes
of Calabi-Yau $3$-folds obtained
from $3$-dimensional hypersurfaces in $4$-dimensional weighted
projective spaces
${\bf P}(\omega_0, \ldots, \omega_4)$. The property $(\ref{main.relation})$
for the construction of B. Greene and R. Plesser \cite{gree} was proved by
S.-S. Roan \cite{roan1}.
\medskip

In the paper  of P. Candelas, X.C. de la Ossa, P.S. Green
and L. Parkes \cite{cand2} the Mirror Symmetry  was applied to
give predictions for the number of rational curves of various
degrees on  general   quintic $3$-folds in ${\bf P}_4$. For degrees
$\leq 3$ these predictions
were confirmed by algebraic geometers  \cite{eil.st,s.kat}.

In \cite{mor1} Morrison has presented
a mathematical review of the calculation of
P. Candelas et al.  \cite{cand2}. Applying analogous method based on
a consideration of the Picard-Fuchs equation, he has found in  \cite{mor2}
similar predictions for the number of rational curves on general  members of
another families of Calabi-Yau $3$-folds with $h^{1,1} = 1$ constructed as
hypersurfaces in weighted projective spaces. Some verifications of  these
these predictions were obtained by S. Katz in \cite{katz.ver}.  The method of
P. Candelas et al.  was also  applied to Calabi-Yau complete intersections
in projective spaces by A. Libgober and J. Teitelbaum \cite{libgober} whose
calculation gave correct predictions for the number of lines and conics.
Analogous results were obtained by physicists  A. Font \cite{font},
A. Klemm and S.Theisen \cite{klemm2,klemm3}.
\bigskip

In this paper we consider families ${\cal F}(\Delta)$ of
Calabi-Yau hypersurfaces  which are compactifications
in $n$-dimensional projective toric
varieties ${\bf P}_{\Delta}$ of smooth affine hypersurfaces
whose equations have a fixed Newton polyhedron $\Delta$ and
sufficiently general coefficients.
Our  purpose  is to describe a general  method for constructing of
candidates for mirrors of Calabi-Yau  hypersurfaces
in toric varieties.
\medskip

It turns out that all results on Calabi-Yau
hypersurfaces and
complete intersections in weighted and projective spaces, in particular
the constructions of mirrors and the computations of
predictions for numbers rational curves, can be extended to
the case of Calabi-Yau hypersurfaces and complete intersections in
toric varieties \cite{batyrev.var,batyrev.quant,bat.cox,bat.straten}.
We remark that the toric technique for resolving singularities
and computing Hodge numbers of  hypersurfaces and
complete intersections was first developed by A.G. Khovansky
\cite{hov.genus,hov}. For the case of  $3$-dimensional varieties with trivial
canonical class, toric methods for resolving quotient singularities
were first applied  by D. Markushevich \cite{mark},  S.-S. Roan
and S.-T. Yau \cite{roan0}.

\bigskip

Let us give an outline of the paper.
\medskip

Section  2 is devoted to basic terminology and  well-known
results on toric varieties.  In this section we fix our notation for the
rest of the paper. We use  two definition of toric varieties.
These definitions  correspond to two kind of combinatorial data,
contravariant and covariant ones.
\medskip

In section 3 we consider general properties of families ${\cal F}(\Delta)$
of hypersurfaces in a toric variety ${\bf P}_{\Delta}$
 satisfying  some regularity
conditions which we call $\Delta$-regularity. These  conditions
imply that the singularities of hypersurfeace are induced only by
singularities of the embient toric variety ${\bf P}_{\Delta}$.
The main consequence of this
fact is that a resolution of singularities
of ${\bf P}_{\Delta}$ immediatelly  gives  rise to a resolution of
singularities of {\em all} $\Delta$-regular hypersurfaces, i.e., we
obtain a {\em simultanious resolution}  of all members of the family
${\cal F}(\Delta)$.
\medskip

In section 4 we investigate lattice  polyhedra  $\Delta$ which
give rise to families of
Calabi-Yau hypersurfaces in ${\bf P}_{\Delta}$. Such a polyhedra
admit a combinatorial characterisation and will be called  {\em reflexive
polyhedra}. There exists the following crusial observation:
\bigskip

{\em If $\Delta \subset M_{\bf Q}$ is a  reflexive
polyhedron, then the corresponding dual polyhedron $\Delta^*$
in dual space $N_{\bf Q}$
is also reflexive.}
\bigskip

Therefore the set of all $n$-dimensional reflexive polyhedra admits
the involution $ \Delta \rightarrow \Delta^*$  which induces the
involution between families of Calabi-Yau varieties:
\[ {MIR}\, :\, {\cal F}({\Delta}) \rightarrow {\cal F}({\Delta}^*). \]
We  use the notaion  $MIR$  for the involution acting on
families of Calabi-Yau hypersurfaces corresponding
to reflexive polyhedra in order to stress the following main
conjecture:
\medskip

{\em The combinatorial involution ${\rm Mir} :
\Delta \rightarrow \Delta^*$ acting on the set
of all reflexive polyhedra of dimension $d$ agrees
with the mirror involution
on conformal field theories
associated to
Calabi-Yau varieties from ${\cal F}({\Delta})$ and  ${\cal F}({\Delta}^*)$. }
\medskip

The next purpose  of the paper is to give  arguments
showing that the  Calabi-Yau families ${\cal F}(\Delta)$  and
${\cal F}(\Delta^*)$ are good candidates to be  mirror symmetric, i.e.,
the  involution
$MIR : {\cal F}({\Delta}) \rightarrow {\cal F}({\Delta}^*)$
agrees  with  properties of the mirror duality in physics.
\medskip

Let ${\hat Z}$ denotes  a {\em maximal projective crepant partial}
de\-sin\-gu\-la\-ri\-za\-tion ({\em MPCP}-de\-sin\-gu\-la\-ri\-za\-tion) of a
projective
Calabi-Yau hypersurface ${\overline Z}$ in ${\bf P}_{\Delta}$.
First, using results of Danilov and Khovanski\^i,  we show  that for
$n \geq 4$ the Hodge number $h^{n-2,1}$  of {\em MPCP}-desingularizations
of Calabi-Yau hypersurfaces  in  the family ${\cal F}(\Delta)$ equals
the Picard number  $h^{1,1}$ of  {\em MPCP}-desingularizations
of  Calabi-Yau hypersurfaces  in the family ${\cal F}(\Delta^*)$.
As a corollary, we obtain the relation (\ref{main.relation}) predicted by
physicists for mirror symmetric Calabi-Yau $3$-folds.
Then, we prove  that
the nonsingular part ${\hat {Z}}$ of $\overline{Z}$ consisting of the
union of $Z_f$ with all $(n-2)$-dimensional affine strata corresponding
$(n-1)$-dimensional edges of $\Delta$ has {\em always} the Euler
characteristic
zero. We use this fact in the proof of a simple combinatorial formula for
the Euler characteristic of Calabi-Yau $3$-folds in terms of
geometry of two $4$-dimensional reflexive polyhedra $\Delta$ and $\Delta^*$.
This formula immediately
implies that for any pair of dual $4$-dimensional
reflexive polyhedra $\Delta$ and $\Delta^*$, Calabi-Yau 3-folds obtained
as {\em MPCP}-desingularizations of Calabi-Yau hypersurfaces
in ${\cal F}(\Delta)$ and in ${\cal F}(\Delta^*)$ have opposite Euler
characteristics.
\medskip

Section  5 is devoted to relations between our method and other
already known methods of explicit constructions of mirrors.

First, we calculate the  one-parameter  mirror family for the family
of $(n-1)$-dimensional hypersurfaces of degree $n+1$ in
$n$-dimensional projective space and show that our result for
quintic therefolds coincides the already known  construction of physicists.

Next,  we investigate the category ${\cal C}_n$
of reflexive pairs $(\Delta, M)$,
where $M$ is an integral lattice of rank $n$ and $\Delta$ is an
$n$-dimensional  reflexive
polyhedron with vertices in $M$. Morphisms in the category ${\cal C}_n$
give rise to relations between families of Calabi-Yau hypersurfaces in
toric varieties.  Namely, the existence of a morphism from
$(\Delta_1, M_1)$ to $(\Delta_2, M_2)$ implies that
Calabi-Yau hypersurfaces in ${\cal F}(\Delta_1)$ consist of
quotients by action of a finite abelian group of Calabi-Yau hypersurfaces
in ${\cal F}(\Delta_2)$.

We prove that if a reflexive polyhedron
$\Delta$ is an $n$-dimensional simplex, then the family
${\cal F}(\Delta)$  consist
of deformations of quotiens of Fermat-type hypersurfaces in
weighted projective spaces.
As a result, we obtain that our method
for  constructions  of mirror candidates  coincides with the method
of Greene and Plesser \cite{gree}, so that we obtain a generalization of
the theorem of S.-S. Roan in \cite{roan1}.
\bigskip

{\bf Acknowledgements.}
The ideas presented in the paper arose  in some preliminary form during my
visiting Japan in  April-June 1991 supported by  Japan Association
for Mathematical Sciences and T\^ohoku University.
I am grateful to  P.M.H. Wilson whose lectures
on Calabi-Yau $3$-folds in Sendai and Tokyo introduced me to this topic.
Geometrical and arithmetical aspects of the theory of
higher-dimensional Calabi-Yau varieties turned out to have deep
connections with my  research interests due to influence of my teachers
V.A. Iskovskih and  Yu. I. Manin.

It is a pleasure to acknowledge helpful recommendations,  and
remarks  concerning  preliminary versions of
this paper from mathematicians D. Dais,  B. Hunt, Y. Kawamata,
D. Morrison, T. Oda, S.-S. Roan, D. van Straten, J. Stienstra,
and  physicists P. Berglund, Ph. Candelas, A. Klemm, J. Louis,
R. Schimmrigk, S. Theisen.

I am very grateful to D. Morrison, B. Greene and P. Aspinwall who
found a serious error in my earlier formulas for the Hodge number
$h^{1,1}$ while working on their  paper \cite{aspin4}.
Correcting these errors, I found an easy proof
of the formula for $h^{n-2,1}$.

I would like to thank the DFG  for the support, and the University of
Essen, especially  H. Esnault and E. Viehweg, for providing
ideal conditions for my work.
\bigskip

\section{The geometry of toric varieties}

\hspace*{\parindent}

We follow notations of V. Danilov  in \cite{dan1}.
Almost all statements of this section are contained
in \cite{oda1} and \cite{reid}.

\medskip

\subsection{Two definitions and basic notations}

\hspace*{\parindent}

First we fix notations used in the contravariant definition
of toric varieties ${\bf P}_{\Delta}$ associated to a lattice
polyhedron $\Delta$.
\bigskip

{$M$}
{ abelian group of rank $n$; }

$\overline{M}$
$ = {\bf Z} \oplus M$;

{${\bf T}$}
$ = ({\bf C}^*)^n = \{ X =
( X_1 , \ldots,  X_n ) \in {\bf C}^n \mid
X_1 \cdots X_n \neq 0 \}$
$n$-dimensional algebraic torus
${\bf T}$  over ${\bf C}$;

{$M_{\bf Q}$}
$= M \otimes {\bf Q}$  the ${\bf Q}$-scalar extension
of $M$;

{$\overline{M}_{\bf Q}$}
$ = {\bf Q} \oplus M_{\bf Q}$;

 {$\Delta$}
a convex $n$-dimensional polyhedron in
$M_{\bf Q}$ with integral vertices (i.e., all vertices of $\Delta$ are
elements of $M$);

 {${\rm vol}_M\, \Delta$ (or $d(\Delta)$)}
the volume of the polyhedron
$\Delta$ relative to the integral lattice $M$, we  call it the
{\em degree} of the polyhedron $\Delta$ relative to $M$;

 {$C_{\Delta}$}
$ = 0 \cup \{ (x_0, x_1, \ldots, x_n ) \in
{\bf Q} \oplus {M}_{\bf Q}
\mid  ({x_1}/{x_0}, \ldots, {x_n}/{x_0})  \in \Delta, \; x_0 >0 \}$
  the  $(n+1)$-dimensional convex
cone  supporting $\Delta$;

 {$S_{\Delta}$}
the the graded  subring
of ${\bf C} \lbrack X_0, X_1^{\pm 1}, \ldots , X_n^{\pm 1} \rbrack$ with
the ${\bf C}$-basis consisting of monomials $X_0^{m_0}X_1^{m_1}
\cdots X_n^{m_n}$
such that $(m_0, m_1, \ldots, m_n) \in C_{\Delta}$;

 {${\bf P}_{\Delta,M}$ (or ${\bf P}_{\Delta}$)}
$ = {\rm Proj}  \, S_{\Delta}$
be an  $n$-dimensional
projective toric  variety corresponding to the graded
ring $S_{\Delta}$;

 {${\cal O}_{\Delta}(1)$}
the ample invertible sheaf on
${\bf P}_{\Delta, M}$ corresponding to the graded $S_{\Delta}$-module
$S_{\Delta}(-1)$;

 {$\Theta$}
be an arbitrary $l$-dimensional polyhedral face of $\Delta$
$(l \leq  n)$;

 {${\bf T}_{\Theta}$}
the corresponding $l$-dimensional ${\bf T}$-orbit in ${\bf P}_{\Delta}$;

 {${\bf P}_{\Theta}$}
the $l$-dimensional toric variety which is the closure of ${\bf T}_{\Theta}$;

 {${\bf P}_{\Delta}^{(i)}$}
$= \bigcup_{{\rm codim}\, \Theta = i}
{\bf P}_{\Theta} = \bigcup_{{\rm codim}\, \Theta \leq i}
{\bf T}_{\Theta }$.
\bigskip

Althought, the definition of toric varieties ${\bf P}_{\Delta}$ associated
with integral polyhedra $\Delta$ is very simple, it is not always convenient.
Notice that a polyhedron $\Delta$  defines not only the corresponding toric
variety ${\bf P}_{\Delta}$, but also the  choice of an
ample invertible sheaf ${\cal O}_{\Delta}(1)$
on ${\bf P}_{\Delta}$ and an  embedding of ${\bf P}_{\Delta}$
into a  projective space. In general,
there exist infinitely many different ample sheaves on ${\bf P}_{\Delta}$.
As a result, there are infinitely many different integral  polyhedra defining
isomorphic toric varieties (one can take, for example,  multiples of
$\Delta$ :  $k \Delta$, $k = 1,2, \ldots )$. If we want to get a
one-to-one
correspondence  between toric varieties and some combinatorial data,
we need covariant definition of toric varieties in terms of fans $\Sigma$
of rational
polyhedral cones.
This  approach to  toric varieties  gives more possibilities,
it allows, for example, to construct affine
and quasi-projective toric varieties as well as  complete
toric varieties which are not quasi-projective.
\bigskip

{$N$}
$= {\rm Hom}\, (M, {\bf Z})$
the dual to $M$ lattice;

 {$\langle *, * \rangle$}\, : \,
$M_{\bf Q} \times N_{\bf Q} \rightarrow
{\bf Q} $ the nondegenerate pairing between the $n$-dimensional
${\bf Q}$-spaces $M_{\bf Q}$ and $N_{\bf Q}$;

 {$\sigma$}
an $r$-dimensional $(0 \leq r \leq n)$ convex rational
polyhedral cone  in $N_{\bf Q}$ having $0 \in N_{\bf Q}$ as vertex;

 {$\check {\sigma}$}
the dual to $\sigma$ $n$-dimensional cone in $M_{\bf Q}$;

 {${\bf A}_{\sigma, N}$ (or ${\bf A}_{\sigma}$)}
$ = {\rm Spec}\, \lbrack \check {\sigma} \cap M \rbrack$ an
$n$-dimensional affine toric variety
associated with the $r$-dimensional cone $\sigma$;

 {$N(\sigma)$}
minimal $r$-dimensional sublattice of $N$ containing $\sigma \cap N$;

 {${\bf A}_{\sigma,N(\sigma)}$}
the $r$-dimensional affine toric variety corresponding to
$\sigma \subset N(\sigma)_{\bf Q}$
(${\bf A}_{N(\sigma)} =
{\bf C}^*)^{n-r} \times {\bf A}_{\sigma, N(\sigma)}$);

 {$\Sigma$}
a finite rational polyhedral fan of cones in $N_{\bf Q}$;

 {$\Sigma^{(i)}$}
the set of all $i$-dimensional cones in $\Sigma$;

 {$\Sigma^{[i]}$}
the subfan of $\Sigma$
consisting of all cones $\sigma \in \Sigma$ such that ${\rm dim}\, \sigma
\leq i $;

 {${\bf P}_{\Sigma,N}$ (or ${\bf P}_{\Sigma}$)}
$n$-dimensional toric variety obtained by glueing of
affine toric varieties ${\bf A}_{\sigma,N}$ where $\sigma \in \Sigma$;

 {${\bf P}^{[i]}_{\Sigma}$}
the open toric subvariety in ${\bf P}_{\Sigma}$
corresponding to $\Sigma^{[i]}$.
\bigskip

Inspite of the fact that the definition of toric varieties via rational
polyhedral fans is more general than via integral polyhedra, we will use both
definitions. The choice of a definition in the
sequel will depend on  questions  we are interested in.
In one situation, it will be more convenient to describe properties
of toric varieties and their subvarieties in terms of integral polyhedra. In
another situation, it will be more convenient to use the language of
rational polyhedral fans.
So it is important to know how one can construct a fan $\Sigma(\Delta)$ from
an integral polyhedron $\Delta$, and how one can construct an integral
polyhedron $\Delta(\Sigma)$ from a rational polyhedral fan $\Sigma$.

The first way $\Delta \Rightarrow \Sigma(\Delta)$ is rather simple:

\begin{prop}
For every  $l$-dimensional face  $\Theta \subset \Delta$, define the convex
$n$-dimensional cone ${\check {\sigma}}(\Theta) \subset M_{\bf Q}$
 consisting of all vectors
$\lambda (p - p')$,  where $\lambda \in{\bf Q}_{\geq 0}$,
$p \in \Delta$, $p' \in \Theta$. Let $\sigma( \Theta) \subset N_{\bf Q}$
be the $(n-l)$-dimensional dual cone
relative to  ${\check {\sigma}}(\Theta) \subset M_{\bf Q}$.

 The set $\Sigma(\Delta)$ of all  cones $\sigma(\Theta)$, where
$\Theta$ runs over all faces of $\Delta$, determines the complete
rational polyhedral fan defining the toric variety ${\bf P}_{\Delta}$.
\label{def.fan}
\end{prop}
\bigskip

The decomposition  of
${\bf P}_{\Delta}$ into a disjoint union of
${\bf T}$-orbits ${\bf T}_{\Theta}$ can be reformulated via cones
$\sigma(\Theta)$
in the fan $\Sigma(\Delta) = \Sigma $ as follows.

\begin{prop}
Let $\Delta$ be an $n$-dimensional $M$-integral polyhedron in $M_{\bf Q}$,
$\Sigma = \Sigma(\Delta)$  the corresponding
complete rational  polyhedral fan
in $N_{\bf Q}$.   Then

{\rm (i)} For any face $\Theta \subset \Delta$,
the affine toric variety ${\bf A}_{{\sigma}(\Theta),N}$ is the
minimal ${\bf T}$-invariant
affine open subset in ${\bf P}_{\Sigma}$ containing  ${\bf T}$-orbit
${\bf T}_{\Theta}$.

{\rm (ii)} Put ${\bf T}_{\sigma}:= {\bf T}_{\sigma(\Theta)}$.
There exists a one-to-one correspondence between
$s$-dimensional cones $\sigma \in \Sigma$ and $(n-s)$-dimensional
${\bf T}$-orbits ${\bf T}_{\sigma}$ such that ${\bf T}_{\sigma'}$ is
contained in the closure of ${\bf T}_{\sigma}$ if and only if
$\sigma$ is a face of $\sigma'$.

{\rm (iii)}
\[ {\bf P}_{\Sigma}^{[i]} = \bigcup_{{\rm dim}\,  \sigma \leq i}
{\bf T}_{\sigma}  \]
is an open ${\bf T}$-invariant subvariety
in ${\bf P}_{\Sigma} = {\bf P}_{\Delta}$, and
\[ {\bf P}_{\Sigma} \setminus {\bf P}_{\Sigma}^{[i]} =
{\bf P}^{(i)}_{\Delta}. \]
\label{orbits.fan}
\end{prop}
\bigskip

Let us now consider another direction $\Sigma \Rightarrow \Delta(\Sigma)$.
In this case, in order to construct $\Delta(\Sigma)$
it is not sufficient to know only a complete fan $\Sigma$.
We need a strictly upper convex support function $h\; :\;  N_{\bf Q}
\rightarrow {\bf Q}$.

\begin{opr}
{\rm Let $\Sigma$ be a rational polyhedral fan and let
$h\; :\;  N_{\bf Q} \rightarrow {\bf Q}$ be a function such that
$h$ is linear on any cone $\sigma \subset \Sigma$. In this situation
$h$ is called a {\em support function} for the fan $\Sigma$.
We call $h$ {\em integral} if $h(N) \subset {\bf Z}$.  We call $h$
{\em upper convex} if  $h(x + x') \leq h(x) + h(x')$ for  all
$x, x' \in N_{\bf Q}$. Finally,  $h$ is called {\em strictly upper convex} if
$h$ is upper convex and
for any two distinct $n$-dimensional cones $\sigma$ and $\sigma'$ in
$\Sigma$ the restrictions $h_{\sigma}$ and $h_{\sigma'}$ of $h$ on $\sigma$
and $\sigma'$ are different  linear functions. }
\end{opr}

\begin{rem}
{\rm By general theory of toric varieties \cite{dan1,oda1}, support functions
one-to-one correspond to  ${\bf T}$-invariant
${\bf Q}$-Cartier divisors $D_h$
on  ${\bf P}_{\Sigma}$. A divisor $D_h$ is a ${\bf T}$-invariant
Cartier divisor (or ${\bf T}$-linearized invertible sheaf) if and only if
$h$ is integral. $D_h$ is numerically effective if and only if $h$ is upper
convex.  Strictly upper convex support functions $h$ on a fan $\Sigma$
one-to-one correspond to ${\bf T}$-linearized {\em ample} invertible
sheaves ${\cal O}(D_h)$  over ${\bf P}_{\Sigma}$.}
\end{rem}

The next proposition describes the construction of a polyhedron $\Delta$  from
a fan $\Sigma$ supplied with a strictly convex support function $h$.

\begin{prop}
Let the convex polyhedron $\Delta = \Delta(\Sigma, h)$ to be defined as follows
\[ \Delta(\Sigma,h) = \bigcap_{\sigma \in \Sigma^{(n)}} ( -h_{\sigma} +
{\check \sigma} ), \]
where integral linear functions $h_{\sigma}: N \rightarrow {\bf Z}$
 are considered as elements of
the lattice $M$. Then one has ${\bf P}_{\Delta} \cong {\bf P}_{\Sigma}$ and
${\cal O}_{\Delta}(1) \cong {\cal O}(D_h)$.
\end{prop}

We have seen already that the construction of an integral polyhedron $\Delta$
from a fan $\Sigma$ such that ${\bf P}_{\Delta} \cong {\bf P}_{\Sigma}$ is
not unique and depends on the choice of an integral
 strictly upper convex support function $h$.
However, there exists  an important case when we can make
a natural choice of $h$.

Indeed, for every
${\bf Q}$-Gorenstein toric variety ${\bf P}_{\Sigma}$ we have the  unique
support function $h_K$ corresponding to the anticanonical
divisor $-K_{\Sigma} = {\bf P}_{\Sigma} \setminus {\bf T}$.

\begin{opr}
{\rm A toric variety ${\bf P}_{\Sigma}$ is said to be a
toric {\em ${\bf  Q}$-Fano variety} if  the anticanonical support function
$h_K$ is strictly upper convex. A toric ${\bf  Q}$-Fano variety is called
a {\em Gorenstein toric Fano variety} if $h_K$ is integral. }
\label{def.fano1}
\end{opr}

With a toric ${\bf Q}$-Fano variety, one can associate two convex
polyhedra:
\[ \Delta({\Sigma}, h_K) = \bigcap_{\sigma \in \Sigma^{(n)}}
( -k_{\sigma} + {\check \sigma} ) \]
and
\[ \Delta^*(\Sigma, h_K) = \{ y \in N_{\bf Q} \mid h_K(y) \leq 1 \}. \]
The polyhedra $\Delta({\Sigma}, h_K)$ and $\Delta^*(\Sigma, h_K)$ belong
to $M_{\bf Q}$ and $N_{\bf Q}$ respectively.

\subsection{Singularities and morphisms of toric varieties}

\hspace*{\parindent}

\begin{opr}
{\rm Let $\sigma \in Sigma$ be an $r$-dimensional cone in $N_{\bf Q}$.
We denote by
$p_{\sigma}$ the unique ${\bf T}/{\bf T}_{\sigma}$-invariant
point on the $r$-dimensional affine toric variety
${\bf A}_{\sigma,N(\sigma)}$ }
\end{opr}

Since  ${\bf A}_{\sigma,N}$ splits into the product
\[ ({\bf C}^*)^{n-s} \times {\bf A}_{\sigma, N(\sigma)}, \]
open analytical neighbourhoods of any two points
$p, p' \in {\bf T}_{\sigma}$ are locally isomorphic. Therefore, it suffices
to investigate the structure of toric singularities at points
$p_{\sigma}$ where $\sigma$ runs over cones of $\Sigma$.

The following two propositions are due to  M. Reid
\cite{reid}.

\begin{prop}
Let $n_1, \ldots, n_r \in N$ $(r \geq s)$ be primitive $N$-integral
generators of all $1$-dimensional faces of an $s$-dimensional
cone $\sigma$.

{\rm (i)} the point $p_{\sigma} \in {\bf A}_{\sigma, N(\sigma)}$ is
${\bf Q}$-factorial
$($or quasi-smooth$)$ if ond only if the cone $\sigma$ is simplicial,
i.e., $r = s$;

{\rm (ii)} the point $p_{\sigma} \in {\bf A}_{\sigma, N(\sigma)}$ is
${\bf Q}$-Gorenstein
if ond only if the elements $n_1, \ldots, n_r$ are contained in
an affine hyperplane
\[ H_{\sigma} \; : \; \{ y \in N_{\bf Q} \mid \langle
k_{\sigma}, y \rangle = 1 \}, \]
for some $k_{\sigma} \in M_{\bf Q}$ $($note that when
${\rm dim}\, \sigma = {\rm dim}\, N$, the element
$k_{\sigma}$ is unique if it exists$)$. Moreover,
${\bf A}_{\sigma, N(\sigma)}$ is Gorenstein if and only if
$k_{\sigma} \in M$.
\label{sing1}
\end{prop}

\begin{rem}
{\rm If ${\bf P}_{\Sigma}$ is a ${\bf Q}$-Gorenstein toric variety, the
elements $k_{\sigma} \in M_{\bf Q}$ ($\sigma \in \Sigma^{(n)}$) define
together
the support function $h_K$ on $N_{\bf Q}$ such that the restriction of
$h_K$ on $\sigma \in \Sigma^{(n)}$ is $k_{\sigma}$. The support function
$h_K$ corresponds to the anticanonical divisor on ${\bf P}_{\Sigma}$.}
\label{anti.des}
\end{rem}

\begin{prop}
Assume that ${\bf A}_{\sigma, N(\sigma)}$ is ${\bf Q}$-Gorenstein {\rm
($see$  \ref{sing1}(ii))}, then

{\rm (i)} ${\bf A}_{\sigma, N(\sigma)}$ has at the point $p_{\sigma}$
at most { terminal} singularity
if and only if
\[ N \cap \sigma \cap \{y \in N_{\bf Q} \mid \langle  k_{\sigma},
y \rangle \leq  1 \}
= \{ 0, n_1, \ldots, n_r \} ; \]

{\rm (ii)} ${\bf A}_{\sigma, N(\sigma)}$ has at the point $p_{\sigma}$ at most
{ canonical } singularity
if and only if
\[ N \cap \sigma \cap \{y \in N_{\bf Q}
\mid \langle  k_{\sigma} , y \rangle <  1 \}
= \{ 0 \}. \]
\label{sing2}
\end{prop}

Using \ref{sing1} and \ref{sing2}, we obtain:

\begin{coro}
Any Gorenstein toric singularity is canonical.
\label{gor.can}
\end{coro}

\begin{opr}
{\rm Let $S$ be a $k$-dimensional simplex in ${\bf Q}^n$ $(k \leq n)$ with
vertices in ${\bf Z}^n$,
$A(S)$ the minimal $k$-dimensional affine ${\bf Q}$-subspace containing
$S$. Denote by ${\bf Z}(S)$ the $k$-dimensional lattice
$A(S) \cap {\bf Z}^n$. We call $P$ {\em elementary} if
$S \cap {\bf Z}(S)$ contains only vertices of $S$. We call $S$
{\em regular}  if the degree of $S$ relative to ${\bf Z}(S)$ is $1$.}
\end{opr}

It is clear that every regular simplex is elementary. The converse is not true
in general. However, there exists the following easy lemma.

\begin{lem}
Every elementary simplex of dimension $\leq 2$ is
regular.
\label{el.reg}
\end{lem}

By \ref{sing1} and  \ref{sing2}, we obtain:

\begin{prop}
Let  ${\bf P}_{\Sigma}$  be a toric variety with only ${\bf Q}$-Gorenstein
singularities, i.e., for any cone $\sigma \in \Sigma$ let the corresponding
element $k_{\sigma} \in M_{\bf Q}$ be well-defined {\rm (\ref{sing1}(ii))}.
Then

{\rm (i)}  ${\bf P}_{\Sigma}$ has only ${\bf Q}$-factorial terminal
singularities if and only if for every cone $\sigma \in \Sigma$ the
polyhedron
\[ P_{\sigma} = \sigma \cap \{ y \in N_{\bf Q} \mid \langle
k_{\sigma}, y \rangle \leq  1 \} \]
is an elementary simplex.

{\rm (ii)} ${\bf P}_{\Sigma}$ is smooth if and only if for every cone
$\sigma \in \Sigma$ the
polyhedron
\[ P_{\sigma} = \sigma \cap \{ y \in N_{\bf Q} \mid \langle
k_{\sigma}, y \rangle \leq  1 \} \]
is a regular simplex.
\label{simp}
\end{prop}

\begin{theo}
Let ${\bf P}_{\Sigma}$ be a toric variety.

{\rm (i)} ${\bf P}_{\Sigma}$ is smooth.

{\rm (ii)} If ${\bf P}_{\Sigma}$ has only
terminal singularities, then the open toric subvariety
${\bf P}_{\Sigma}^{[2]}$ is smooth.

{\rm (iii)} If ${\bf P}_{\Sigma}$ has only
Gorenstein ${\bf Q}$-factorial
terminal singularities, then the open toric subvariety
${\bf P}_{\Sigma}^{[3]}$ is smooth.
\label{codim2,3}
\end{theo}

We recall a combinatorial characterization of toric  morphisms
between toric varieties.
\bigskip

Let $\phi \; : \; N' \rightarrow N$ be a morphism of lattices, $\Sigma$ a
fan in $N_{\bf Q}$,  $\Sigma'$ a fan in $N'_{\bf Q}$. Suppose that for each
$\sigma' \in \Sigma'$ we can find a $\sigma \in \Sigma$ such that
$\varphi(\sigma') \subset \sigma$. In this situation there arises  a morphism
of toric varieties
\[ \tilde {\phi} \; : \; {\bf P}_{\Sigma',N'} \rightarrow
{\bf P}_{\Sigma, N}. \]

\begin{exam} {Proper birational morphisms.} {\rm Assume that
$\varphi$ is an isomorphism of lattices and $\phi(\Sigma')$
is a subdivision of $\Sigma$, i.e., every cone $\sigma \in \Sigma$ is a
union of cones of $\varphi(\Sigma')$. In this case  $\tilde {\phi}$ is
a proper birational  morphism. Such a morphism we will use  for
constructions of desingularizations of toric singularities. }
\label{biration}
\end{exam}
\medskip

\begin{opr}
{\rm Let $\varphi:  W' \rightarrow W$ be a proper birational
morphism of normal ${\bf Q}$-Gorenstein algebraic varieties. The morphism
$\varphi$ is called {\em crepant} if $\varphi^* K_{W} = K_{W'}$
($K_W$ and $K_{W'}$ are canonical divisors  on $W$ and $W'$ respectively). }
\end{opr}
\medskip

Using the description in \ref{anti.des}
of support functions corresponding to anticanonical
divisors on ${\bf P}_{\Sigma}$ and ${\bf P}_{\Sigma'}$, we obtain.

\begin{prop}
A proper birational  morphism of ${\bf Q}$-Gorenstein toric
varieties
\[ \tilde {\phi} \; : \; {\bf P}_{\Sigma'} \rightarrow
{\bf P}_{\Sigma} \]
is crepant if and only if for every cone $\sigma \in \Sigma^{(n)}$ all
$1$-dimensional cones $\sigma'\in  \Sigma'$ which are
contained in $\sigma$ are
generated by primitive integral elements from $N \cap H_{\sigma}$
{\rm ($see$ \ref{sing1}(ii))}.
\label{crit.crep}
\end{prop}

\begin{opr}
{\rm Let $\varphi\; : \; W' \rightarrow W$ be a projective birational
morphism of normal ${\bf Q}$-Gorenstein algebraic varieties. The morphism
$\varphi$ is called a {\em maximal projective crepant partial
desingularization} ({\em MPCP-desingularization}) of $W$
if $\varphi$ is crepant
and $W'$ has only ${\bf Q}$-factorial terminal singularities.}
\end{opr}
\medskip

Our next purpose is to define some  combinatorial notions which we
use to construct $MPCP$-desingularizations of
Gorenstein toric varieties (see \ref{crep.fano}).
\medskip

\begin{opr}
{\rm Let $A$ be a finite subset in $\Delta \cap {\bf Z}^n$. We call $A$
{\em admissible} if it  contains
all vertices of the integral polyhedron $\Delta$. }
\end{opr}

\begin{opr}
{\rm Let $A$ be an admissible subset in $\Delta \cap {\bf Z}^n$.
By an  {\em A-triangulation} of $\Delta \subset
{\bf Q}^n$ we mean a finite collection ${\cal T} = \{ \theta \}$ of simplices
with vertices in $A$ having the following properties:

(i) if $\theta'$ is a face of a simplex $\theta \in {\cal T}$, then
$\theta' \in {\cal T}$;

(ii) the vertices of each simplex $\theta \in {\cal T}$ lie in $\Delta \cap
{\bf Z}^n$;

(iii) the intersection of any two simplices $\theta_1, \theta_2 \in
{\cal T}$ either is empty or is a common face of both;

(iv) $\Delta = \bigcup_{\theta \in {\cal T}} \theta$;

(v) every element of  $A$ is a vertex of some simplex
$\theta \in {\cal T}$. }
\label{def.triang}
\end{opr}

\begin{opr}
{\rm An $A$-triangulation ${\cal T}$ of an integral convex polyhedron
$\Delta \subset
{\bf Q}^n$ is called {\em maximal} if $A = \Delta \cap {\bf Z}^n$.   }
\end{opr}

\begin{rem}
{\rm Note that ${\cal T}$ is a maximal triangulation of $\Delta$
if and only if every simplex $\theta \in {\cal T}$ is elementary. Therefore,
if ${\cal T}$ is a maximal triangulation of $\Delta$, then
for any face $\Theta \subset \Delta$ the induced triangulation of
$\Theta$ is also maximal.}
\label{rem.max}
\end{rem}
\medskip

Assume that  $A$ is admissible.
 Denote by ${\bf Q}^A$  the
${\bf Q}$-space  of functions from $A$ to ${\bf Q}$. Let ${\cal T}$ be a
triangulation of $\Delta$. Every element $\alpha \in {\bf Q}^A$ can be
uniquely extended to a piecewise linear function
$\alpha({\cal T})$ on $\Delta$
such that the restriction of $\alpha({\cal T})$ on every simplex $\theta \in
{\cal T}$ is an affine linear function.

\begin{opr}
{\rm Denote by $C({\cal T})$ the convex cone in ${\bf Q}^A$ consisting of
elements $\alpha$ such that $\alpha({\cal T})$ is an upper
convex function. We say that a
triangulation ${\cal T}$ of $\Delta$ is {\em projective} if the cone
$C({\cal T})$ has a nonempty  interior. In other words, ${\cal T}$ is
projective if and only if there exists a strictly upper convex function
$\alpha({\cal T})$. }
\label{opr.t.proj}
\end{opr}

\begin{prop}
{\rm \cite{gelf}}
Let $\Delta$ be an integral polyhedron. Take an admissible subset $A$ in
$\Delta \cap {\bf Z}^n$. Then $\Delta$ admits at least one projective
$A$-triangulation, in particular, there exists at least one maximal
projective triangulation of $\Delta$.
\label{triang.ex}
\end{prop}

\begin{rem}
{\rm In the paper  Gelfand, Kapranov, and Zelevinsky \cite{gelf}, it has
been used  the notion of a {\em regular triangulation of a polyhedron}
which is called in this paper
a {\em projective triangulation}. The reason for such a change of
the terminology we  will see in \ref{crep.fano}.}
\end{rem}

\begin{exam} {\rm
{\em Finite morphisms.} Assume that $\phi$  is injective,
$\phi(N')$ is a sublattice of finite index
in $N$, and $\phi(\Sigma') = \Sigma$. Then  $\tilde {\phi}$ is a finite
surjective  morphism of toric varieties. This  morphism induces an
\^etale covering of open subsets
\[ {\bf P}_{\Sigma'}^{[1]} \rightarrow {\bf P}_{\Sigma}^{[1]}  \]
with the Galois group ${\rm Coker}\,\lbrack N'\rightarrow N
\rbrack$. }
\label{final}
\end{exam}

Recall  the following  description of the fundamental group
of toric varieties.

\begin{prop}
The fundamental group of a toric variety ${\bf P}_{\Sigma}$ is isomorphic
to the quotient of $N$ by sum of all sublattices $N(\sigma)$ where $\sigma$
runs over all cones $\sigma \in \Sigma$. In particular, the fundamental
group of the non-singular open toric subvariety ${\bf P}_{\Sigma}^{[1]}$
is isomorphic to the quotient of $N$ by sublattice spanned by all
primitive integral generators  of $1$-dimensional cones $\sigma \in
\Sigma^{(1)}$.
\label{fund.group}
\end{prop}
\bigskip

We will use the following characterization of toric Fano varieties in terms of
he polyhedra $\Delta({\Sigma}, h_K)$ and
$\Delta^*({\Sigma}, h_K)$.

\begin{prop}
A complete toric variety ${\bf P}_{\Sigma}$ with only ${\bf Q}$-Gorenstein
singularities is a ${\bf Q}$-Fano variety {\rm ($see$ \ref{def.fano1})}
if and only if $\Delta({\Sigma}, h_K)$ is an
$n$-dimensional polyhedron with vertices $-k_{\sigma}$
one-to-one corresponding
to $n$-dimensional cones $\sigma \in \Sigma$.  In this situation,

{\rm (i)} ${\bf P}_{\Sigma}$ is a Fano variety with only Gorenstein
singularities if and only if all vertices of $\Delta({\Sigma}, h_K)$
belong to $M$, in particular,
${\bf P}_{\Sigma} = {\bf P}_{\Delta({\Sigma, h_K})}$.

{\rm (ii)} ${\bf P}_{\Sigma}$ is a smooth Fano variety if and only if for
every $n$-dimensional cone $\sigma \in \Sigma^{(n)}$ the polyhedron
$\Delta^*({\Sigma}, h_K) \cap \sigma$ is a  regular simplex
of dimension $n$.
\label{smooth.fano}
\end{prop}

Now we come to the most important statement which will be used
in the sequel.

\begin{theo}
Let ${\bf P}_{\Sigma}$ be a toric Fano variety with only Gorenstein
singularities. Then ${\bf P}_{\Sigma}$ admits at least one
MPCP-desingulari\-zation
\[ \tilde {\phi} \; : \; {\bf P}_{\Sigma'} \rightarrow
{\bf P}_{\Sigma}. \]
Moreover,  MPCP-desingulari\-zations of
${\bf P}_{\Sigma}$ are defined by maximal projective triangulations of
the polyhedron
$ \Delta^*(\Sigma, h_K)$,
where $h_K$ the integral strictly upper convex support function associated
with the anticanonical divisor ${\bf P}_{\Sigma} \setminus {\bf T}$ on
${\bf P}_{\Sigma}$.
\label{crep.fano}
\end{theo}

\proof  Define  a finite subset $A$ in $N$ as follows
\[ A = \{ y \in N \mid h_K(y) \leq 1 \} = N \cap \Delta^*({\Sigma}, h_K). \]
It is clear that $A$ is an admissible subset of $\Delta^*(\Sigma, h_K)$.
By \ref{triang.ex}, there exists
at  least one projective $A$-triangulation ${\cal T}$ of
$\Delta^*(\Sigma, h_K)$. Let
$B$ be the boundary of $\Delta^*(\Sigma, h_K)$. For every
simplex $\theta \in B$, we construct a
convex cone $\sigma_{\theta}$ supporting $\theta$. By definition
\ref{def.triang}, the set of all cones $\sigma_{\theta}$
($\theta \in {\cal T} \cap B$) defines a fan $\Sigma'$
which is a subdivision of $\Sigma$.

Since generators of $1$-dimensional
cones of $\Sigma$ are exactly elements of $A \cap B$, the morphism
${\bf P}_{\Sigma'} \rightarrow {\bf P}_{\Sigma}$ is crepant
(see \ref{crit.crep}).

By \ref{opr.t.proj}, there exists a strictly upper convex function
$\alpha({\cal T})$. We can also assume that $\alpha({\cal T})$ has
zero value at $0 \in N$. Then $\alpha({\cal T})$ defines a
strictly convex support function for the fan $\Sigma'$. Thus,
${\bf P}_{\Sigma'}$ is projective.

By \ref{simp}(i) and \ref{rem.max}, we obtain that the morphism
${\bf P}_{\Sigma'} \rightarrow {\bf P}_{\Sigma}$ is
a $MPCP$-desingulari\-zation.

By similar arguments, one can prove that any $MPCP$-desingulari\-zation
defines
a maximal projective triangulation of $\Delta^*(\Sigma, h_K)$.
\bigskip

\section{Hypersurfaces in toric varieties}

\subsection{Regularity conditions for hypersurfaces}

\hspace*{\parindent}

A Laurent polynomial $f = f(X)$  is a finite linear
combination of elements of $M$
\[ f(X)  = \sum c_m X^m \]
with complex coefficients $c_m$. The Newton polyhedra $\Delta(f)$
of $f$ is the convex
hull in $M_{\bf Q} = M \otimes {\bf Q}$ of all elements $m$ such
that  $c_m \neq 0$. Every Laurent polynomial $f$ with the Newton
polyhedron $\Delta$ defines the  affine hypersurface
\[ Z_{f, \Delta} = \{ X \in {\bf T}  \mid f(X) = 0 \}. \]
If we work with a fixed Newton polyhedron $\Delta$, we denote
$Z_{f, \Delta}$ simply by $Z_f$.
\bigskip

Let $\overline{Z}_{f,\Delta}$ be the
closure of $Z_{f, \Delta} \subset {\bf T}$ in  ${\bf P}_{\Delta}$.
For any face $\Theta \subset \Delta$, we put
$Z_{f, \Theta} = \overline{Z}_{f, \Delta} \cap {\bf T}_{\Theta}$.
So we obtain  the  induced
decomposition into the disjoint union
\[ \overline{Z}_{f, \Delta} = \bigcup_{\Theta \subset \Delta}
Z_{f, \Theta}. \]
\medskip

\begin{opr}
{\rm
Let $L(\Delta)$ be the space of all Laurent polynomials  with a fixed
Newton polyhedron $\Delta$.
A Laurent polynomial $f \in L(\Delta)$ and the corresponding
hypersurfaces $Z_{f, \Delta} \subset {\bf T}_{\Delta}$,
$\overline{Z}_{f, \Delta} \subset {\bf P}_{\Delta}$ are  said to be
${\Delta}$-{\em regular} if for every  face
$\Theta \subset \Delta$ the affine variety $Z_{f,\Theta}$
is empty  or a smooth subvariety of codimension 1 in ${\bf T}_{\Theta}$.
Affine varieties $Z_{f,\Theta}$  are called  {\em the strata} on
$\overline{Z}_{f, \Delta}$ associated with  faces $\Theta \subset
\Delta$.}
\label{d.nondeg}
\end{opr}

\begin{rem}
{\rm Notice that if $f$ is $\Delta$-regular,
then ${Z}_{f, \Theta} = \emptyset$  if and only if ${\rm dim}\, \Theta = 0$,
i.e., if and only if $\Theta$  is a vertex of $\Delta$. }
\label{zero.orb}
\end{rem}
\bigskip

Since the space $L(\Delta)$ can be identified with the space of global
sections of the ample sheaf ${\cal O}_{\Delta}(1)$ on
${\bf P}_{\Delta}$, using  Bertini theorem, we obtain:

\begin{prop}
The set of $\Delta$-regular hypersurfaces is a Zariski open subset
in $L(\Delta)$.
\end{prop}

We extend the notion of $\Delta$-regular hypersurfaces in
${\bf P}_{\Delta}$ to the case of hypersurfaces in general toric varieties
${\bf P}_{\Sigma}$ associated with rational polyhedral fans $\Sigma$.
\bigskip

\begin{opr}
{\rm Let $\overline{Z}_{f, \Sigma}$ be the closure in ${\bf P}_{\Sigma}$
of an  affine hypersurface $Z_f$ defined by a Laurent polynomial $f$. Consider
the  induced decomposition into the disjoint union
\[ \overline{Z}_{f, \Sigma } = \bigcup_{\sigma  \in  \Sigma}
Z_{f, \sigma},  \]
where ${Z}_{f, \sigma } = \overline{Z}_{f, \Sigma } \cap
{\bf T}_{\sigma}$.
A Laurent polynomial $f$ and the corresponding
hypersurfaces $Z_{f} \subset {\bf T}$,
$\overline{Z}_{f, \Sigma} \subset {\bf P}_{\Sigma}$ are  said to be
${\Sigma}$-{\em regular} if for every  $s$-dimensional  cone
$\sigma  \in  \Sigma$ the variety $Z_{f,\sigma}$
is empty  or  a smooth subvariety of codimension 1 in ${\bf T}_{\sigma}$.
In other words, $\overline{Z}_{f, \Sigma}$
has only transversal intersections with all ${\bf T}$-orbits
${\bf T}_{\sigma}$ ($\sigma \in \Sigma$).
Affine varieties $Z_{f,\sigma}$  are called  {\em strata} on
$\overline{Z}_{f, \Sigma}$ associated with the cones  $\sigma  \subset
\Sigma$.

Denote by $Z_{f, \Sigma}^{[i]}$ the open subvariety
in $\overline{Z}_{f, \Sigma}$ defined as follows
\[ Z_{f, \Sigma}^{[i]} : = \bigcup_{\sigma  \in  \Sigma^{[i]}} Z_{f, \sigma} =
\overline{Z}_{f, \Sigma} \cap {\bf P}_{\Sigma}^{[i]}.   \]}
\end{opr}
\medskip

Let $\sigma$ be an $s$-dimensional cone in $\Sigma$.
If we apply the  implicit function theorem and the standard
criterion of smoothness
to the affine hypersurface $Z_{f,\sigma} \subset {\bf T}_{\sigma}$ contained
in the open $n$-dimensional affine toric variety
\[ {\bf A}_{\sigma, N}   \cong ({\bf C}^*)^{n-s} \times
{\bf A}_{\sigma, N(\sigma)},  \]
then we obtain:
\medskip

\begin{theo}
Small analytical neighbourhoods of points   on  a
$(n -s -1)$-dimensional stratum $Z_{f, \sigma} \subset
\overline{Z}_{f, \Sigma}$ are analytically isomorphic
to products  of a $(s-1)$-dimensional open ball and a small analytical
neighbourhood  of
the point $p_{\sigma}$ on  the $(n -s)$-dimensional affine toric
variety ${\bf A}_{\sigma, N(\sigma)}$.
\label{anal}
\end{theo}

For the case of $\Delta$-regular hypersurfaces in a toric variety
${\bf P}_{\Delta}$ associated with an integral polyhedron $\Delta$,
one gets from \ref{def.fan} the following.

\begin{coro}
For any $l$-dimensional face $\Theta \subset \Delta$, small
analytical neighbourhoods of points  on the
$(l -1)$-dimensional stratum $Z_{f, \Theta} \subset
\overline{Z}_{f, \Delta}$ are analytically isomorphic
to products  of a $(l-1)$-dimensional open ball and a small  analytical
neighbourhood of
the point $p_{\sigma(\Theta)}$ on  the $(n -l)$-dimensional affine toric
variety ${\bf A}_{\sigma(\Theta), N(\sigma(\Theta))}$.
\label{anal2}
\end{coro}

Applying  \ref{codim2,3}, we also conclude:

\begin{coro}
For any $\Sigma$-regular hypersurface $\overline{Z}_{f, \Sigma} \subset
{\bf P}_{\Sigma}$, the
open subset $Z_{f, \Sigma}^{[1]}$ consists of smooth points of
$\overline{Z}_{f, \Sigma}$. Moreover,

{\rm (i)} $Z_{f, \Sigma}^{[2]}$ consists of smooth points if
${\bf P}_{\Sigma}$ has only terminal singularities.

{\rm (ii)} $Z_{f, \Sigma}^{[3]}$ consists of smooth points if
${\bf P}_{\Sigma}$ has only ${\bf Q}$-factorial Gorenstein
terminal singularities.

{\rm (iii)} $Z_{f, \Sigma}^{[n-1]} = \overline{Z}_{f, \Sigma}$
is  smooth  if and only if ${\bf P}_{\Sigma}^{[n-1]}$ is smooth.
\label{sing.hyp}
\end{coro}
\bigskip

\subsection{Birational and finite morphisms of hypersurfaces}

\hspace*{\parindent}

\begin{prop}
Let $\phi\; : \Sigma' \rightarrow \Sigma$  be a subdivision of a fan $\Sigma$,
\[ \tilde {\phi} \; : \; {\bf P}_{\Sigma'} \rightarrow
{\bf P}_{\Sigma} \]
the corresponding proper birational morphism. Then
for any $\Sigma$-regular hypersurface $\overline{Z}_f \subset
{\bf P}_{\Sigma}$ the hypersurface $\overline{Z}_{{\tilde {\phi}}^{*}f}
\subset {\bf P}_{\Sigma'}$  is  $\Sigma'$-regular.
\label{subdiv.hyp}
\end{prop}

\proof  The statement follows from the fact that for any cone $\sigma' \in
\Sigma'$ such that $\phi(\sigma') \subset \sigma \in \Sigma$, one has
\[ Z_{{\tilde {\phi}}^* f, \sigma'} \cong Z_{f, \sigma} \times
({\bf C}^*)^{{\rm dim}\, \sigma - {\rm dim}\, \sigma'}. \]
\bigskip

One can use \ref{sing.hyp} and \ref{subdiv.hyp} in order to construct
partial desingularizations of hypersurfaces $\overline{Z}_{f, \Sigma}$.

\begin{prop}
Let ${\bf P}_{\Sigma}$ be a projective toric variety with only Gorenstein
singularities. Assume that
\[ \tilde {\phi} \; : \; {\bf P}_{\Sigma'} \rightarrow
{\bf P}_{\Sigma} \]
is a MPCP-desingularization of ${\bf P}_{\Sigma}$. Then
$\overline{Z}_{{\tilde \phi}^*f, \Sigma'}$ is  a
MPCP-desingularization of $\overline{Z}_{f, \Sigma}$.
\label{max.crep}
\end{prop}

\proof  By \ref{sing1}, \ref{sing2}, \ref{anal},
$\overline{Z}_{f, \Sigma}$ has at most Gorenstein singularities and
$\overline{Z}_{{\tilde \phi}^*f, \Sigma'}$  has at most ${\bf Q}$-factorial
terminal singularities. It suffuces now to apply the adjunction formula.
\medskip

\begin{prop}
Let $\phi \; : \; N' \rightarrow N$ be a surjective homomorphism
of $n$-dimensional lattices, $\Sigma$ a
fan in $N_{\bf Q}$,  $\Sigma'$ a fan in $N'_{\bf Q}$. Assume that
$\phi(\Sigma') = \Sigma$. Let  \[ \tilde {\phi} \; : \;
{\bf P}_{\Sigma',N'} \rightarrow {\bf P}_{\Sigma, N} \]
be the corresponding finite surjective morphism of toric varieties.
Then for any $\Sigma$-regular hypersurface $\overline{Z}_f$ in
${\bf P}_{\Sigma}$ the hypersurface
\[ {\tilde {\phi}}^{-1} (\overline{Z}_f)  =
\overline{Z}_{{\tilde {\phi}}^{*}f} \]
is $\Sigma'$-regular.
\label{galois}
\end{prop}

\proof  It is sufficient to observe that for any cone $\sigma \in
\Sigma \cong \Sigma'$ the affine variety ${Z}_{{\tilde {\phi}}^{*}f,
\sigma}$ is
an \^etale covering of $Z_{f, \sigma}$ whose Galois group is  isomorphic
to the cokernel of the homomorphism $N'(\sigma) \rightarrow N(\sigma)$
(see  \ref{final}).
\bigskip

\subsection{The Hodge structure of hypersurfaces}

\hspace*{\parindent}

We are interested now in the calculation of the
cohomology groups of $\Delta$-regular
hypersurfaces $\overline{Z}_f$ in toric varieties ${\bf P}_{\Delta}$.
The main difficulty in this calculation is connected with
singularities of $\overline{Z}_f$. So we will try to avoid singularities and
to calculate cohomology groups not only of compact compact hypersurfaces
$\overline{Z}_f$, but also ones of some naturally defined smooth open
subsets in $\overline{Z}_f$. For the last purpose, it is more convenient
to use cohomology with compact supports which we denote by $H^i_c(*)$.
\medskip

First we note that there exists the following
analog of the  Lefschetz theorem for $\Delta$-regular hypersurfaces
proved by Bernstein, Danilov and Khovanski\^i \cite{dan.hov}.

\begin{theo}
For any open toric subvariety $U \subset {\bf P}_{\Delta}$,
the Gysin homomorphism
\[ H^i_c ( \overline{Z}_f \cap U) \rightarrow H^{i+2}_c (U) \]
is bijective for $i > n -1$ and injective for $i = n - 1$.
\label{Lef}
\end{theo}

Using this theorem, one can often reduce the calculation  of
cohomology groups $H^i$ (or $H^i_c$) to the "interesting"
case $i = n-1$. We consider below several typical examples of
such a situation.

\begin{exam}
{\rm Let $U$ be a smooth open toric subvariety in ${\bf P}_{\Delta}$
(e.g., $U = {\bf T}$). Then
$V = U \cap \overline{Z}_f$ is a smooth affine open subset in
$\overline{Z}_f$. By general properties of Stein varieties, one has
$H^{i}(V) = 0$ for $i > n-1$.  Since the calculation of cohomology groups of
smooth affine toric varieties is very simple, we obtain a complete
information about all cohomology groups except for $i = n-1$ using the
following property. }
\end{exam}

\begin{prop}
Let  $W$ be a quasi-smooth irreducible $k$-dimensional algebraic variety. Then
there exists the Poincare pairing
\[ H_c^i (W) \otimes H^{2k -i} (W) \rightarrow H^{2k}_c (W) \cong {\bf C}.  \]
This pairing is compatible with Hodge structures, where $H^{2k}_c(W)$ is
assumed to be a $1$-dimensional ${\bf C}$-space of the Hodge type $(k,k)$.
\label{duality}
\end{prop}

If we take $U = {\bf T}$, then $V = Z_f$. The Euler characteristic  of
$Z_f$ was calculated by
Bernstein, Khovanski\^i and Kushnirenko \cite{kuch}, \cite{hov.genus}.

\begin{theo}
$ e(Z_f) = \sum_{i \geq 0} (-1)^i {\rm dim}\, H^i(Z_f)  =
(-1)^{n-1} d_M (\Delta)$.
\label{euler}
\end{theo}

In particular, we obtain:

\begin{coro}
The dimension of $H^{n-1}(Z_f)$ is equal to $d_M(\Delta) + n -1$.
\end{coro}

\begin{opr}
Let $P$ be a compact convex subset in $M_{\bf Q}$.  Denote by $l(P)$
the number of integral points in $P \cap M$, and by $l^*(P)$
the number of integral points in the interior of $P$.
\end{opr}

There exist general formulas  for the Hodge-Deligne numbers
$h^{p,q}(Z_f)$ of the mixed  Hodge structure in the $(n-1)$-th cohomology
group of an arbitrary
$\Delta$-regular affine hypersurface in ${\bf T}$ (Danilov and Khovanski\^i
\cite{dan.hov}).
For the numbers $h^{n-2,1}(Z_f)$ and $h^{n-2,0}(Z_f)$,  we get the following
(see \cite{dan.hov}, 5.9).

\begin{prop}
Let ${\rm dim}\, \Delta = n \geq 4$, then
\[ h^{n-2,1} (Z_f) + h^{n-2,0} (Z_f) = l^*(2\Delta) - (n+1)l^*(\Delta), \]
\[ h^{n-2,0} (Z_f) =  \sum_{{\rm codim}\,  \Theta =1} l^*(\Theta).\]
\label{hd.dh}
\end{prop}
\bigskip

We will use in the sequel the following properties
of the numbers $h^{p,q}(H^k_c(Z_f))$ of affine hypersurfaces
for cohomology with compact supports (see \cite{dan.hov}):

\begin{prop}
The Hodge-Deligne numbers $h^{p,q}(H^k_c(Z_f)$ of $\Delta$-regular
$(n-1)$-di\-men\-si\-onal affine hypersurfaces $Z_f$ satisfy the properties:

{\rm (i)} $h^{p,q}(H^k_c(Z_f) = 0$ for $p \neq q$ and $k > n-1$;

{\rm (ii)} $h^{p,q}(H^k_c(Z_f)) = 0$ for  $k < n-1$;

{\rm (iii)} $h^{p,q}(H^{n-1}_c(Z_f) = 0$ for  $p+ q  > n-1$.
\label{properties.hodge}
\end{prop}
\bigskip

Although we consider the mixed Hodge structure in the cohomology group
of the affine hypersurface $Z_f$, we  get eventually some information about
the Hodge numbers of compactifications of $Z_f$. From general
properties of the mixed Hodge structures \cite{deligne}, one obtains:

\begin{theo}
Let $\overline{Z}$ be a smooth compactification of a smooth affine
$(n-1)$-dimensional variety $Z$ such that the complementary set
$\overline{Z} \setminus Z$ is a normal crossing divisor. Let
\[ j\;: \;  Z \hookrightarrow \overline{Z} \]
be the corresponding embedding. Denote by
\[ j^* \;:\; H^{n-1} (\overline{Z}) \rightarrow H^{n-1} (Z) \]
the induced mapping of cohomology groups. Then
the weight subspace ${\cal W}_{n-1}H^{n-1}(Z)$ coincides with
the image  $j^*( H^{n-1}(\overline{Z}))$.  In particular, one has
the following inequalities between the Hodge numbers of $\overline{Z}$ and
the Hodge-Deligne numbers of $Z$:
\[ h^{n-1-k,k}(\overline{Z}) \geq h^{n-1-k,k}(H^{n-1}(Z))\;\;
(0 \leq k \leq n-1). \]
\label{inequal.h}
\end{theo}
\bigskip

\section{Calabi-Yau hypersurfaces in toric varieties}
\medskip

\subsection{Reflexive polyhedra and reflexive pairs}

\hspace*{\parindent}

\begin{opr}
{\rm If  $P$ is an arbitrary
compact convex set in $M_{\bf Q}$ containing the
zero vector $0 \in M_{\bf Q}$ in its interior, then we call
\[ P^* = \{ y \in N_{\bf Q} \mid \langle x, y \rangle \geq -1, \;
{\rm for \; all} \; x \in P  \}. \]
{\em the dual set}  relative to $P$.
 }
\label{dual}
\end{opr}
\bigskip

The dual set $P^*$ is  a convex compact subset in $N_{\bf Q}$
containing  the vector zero $0 \in N_{\bf Q}$ in its interior.
Obviously, one has $(P^*)^* = P$.
\bigskip

\begin{exam}
{\rm Let $E$ be a Euclidian $n$-dimensional space, $\langle *, * \rangle$
the corresponding scalar product,
\[ P = \{ {x} \mid \langle { x}, { x} \rangle \leq R \} \]
the  ball of  radius $R$. Using the scalar product, we can identify
the dual space $E^*$ with $E$. Then  the dual set $P^*$ is
the  ball of radius $1/R$. }
\end{exam}

\begin{prop}
Let $P \subset M_{\bf Q}$ be a convex set containing $0$ in its interior,
$C_P \subset \overline{M}_{\bf Q}$ the convex cone supporting $P$, $
\overline{N}_{\bf Q} = {\bf Q} \oplus N_{\bf Q}$ the dual space,
$C_{P^*} \subset \overline{N}_{\bf Q}$
the convex cone supporting $P^* \subset N_{\bf Q}$. Then $C_{P^*}$ is the
dual cone relative to $C_P$.
\label{dual.cones}
\end{prop}

\proof  Let $\overline{x} = (x_0, x) \in \overline{M}_{\bf Q}$,
$\overline{y} = (y_0, y) \in \overline{M}_{\bf Q}$.
Since $x_0$ and $y_0$ are positive, one has
\[ \overline{x} \in C_P \Leftrightarrow {x}/{x_0} \in P \;\; {\rm and }\;\;
\overline{y} \in C_{P^*} \Leftrightarrow {y}/{y_0} \in P^*. \]
Therefore,
$\langle \overline{x}, \overline{y} \rangle = x_0y_0 +
\langle {x}, {y} \rangle \geq 0 $
if and only if
$\langle {x}/{x_0}, {y}/{y_0} \rangle \geq -1$.
\bigskip

\begin{opr}
{\rm Let
$H$ be a rational affine hyperplane in $M_{\bf Q}$, $p \in M_{\bf Q}$ an
arbitrary integral point. Assume that $H$ is affinely generated by
integral points $H \cap M$, i.e., there exists  a primitive integral element
$l \in N$  such that for some integer $c$
\[ H = \{ x \in M_{\bf Q} \mid \langle x, l \rangle = c \}. \]
Then the absolute value
$\mid c - \langle p , l \rangle \mid$ is called the {\em integral distance}
between $H$ and $p$. }
\end{opr}

\begin{opr}
{\rm Let $M$ be an integral $n$-dimensional lattice in $M_{\bf Q}$,
$\Delta$ a convex integral polyhedron in $M_{\bf Q}$ containing the zero
 $0 \in M_{\bf Q}$  in its  interior.
The pair $(\Delta, M)$ is called {\em reflexive} if
the integral distance between
 $0$ and all affine hyperplanes   generated by  $(n -1)$-dimensional faces
of $\Delta$ equals 1.

If $(\Delta, M)$ is a reflexive pair, then we call $\Delta$ a {\em
reflexive polyhedron}.}
\label{inver.p}
\end{opr}
\bigskip

The following simple property of reflexive polyhedra is very important.
\medskip

\begin{theo}
Suppose that $(\Delta, M)$ is a reflexive pair. Then $(\Delta^*, N)$ is
again a reflexive pair.
\end{theo}

\proof  Let $\Theta_1, \ldots , \Theta_k$ be $(n-1)$-dimensional faces of
$\Delta$, $H_1, \ldots , H_k$ the corresponding affine hyperplanes.
By \ref{inver.p}, there exist integral elements
$l_1, \ldots , l_k \in N_{\bf Q}$ such that for all $1 \leq i \leq k$
\[ \Theta_i = \{ x \in \Delta \mid \langle x , l_i \rangle = 1 \},\;\;
 H_i = \{ x \in M_{\bf Q} \mid \langle x , l_i \rangle = 1 \}. \]
Therefore,
\[ \Delta = \{ x \in M_{\bf Q}  \mid \langle x , l_i \rangle \leq 1 \;
(1 \leq i \leq k ) \}. \]
So  $\Delta^*$ is a convex hull of
the integral points $l_1, \ldots , l_k$, i.e., $\Delta^*$ is an
integral polyhedron.

Let $p_1, \ldots , p_m$ be vertices of $\Delta$. By  \ref{dual},
for any $j$ ($1 \leq j  \leq m$)
\[ \Xi_j = \{ y \in \Delta^*  \mid \langle p_j , y \rangle = 1 \} \]
is a $(n-1)$-dimensional face of $\Delta^*$.
Thus, $\Delta^*$ contains $0 \in N_{\bf Q}$  in its interior,
and the integral distance between $0$ and every  $(n-1)$-dimensional
affine linear subspace generated by $(n-1)$-dimensional
faces of $\Delta^*$ equals $1$.
\bigskip

We can establish the following one-to-one correspondence between faces of
the polyhedra $\Delta$ and $\Delta^*$.

\begin{prop}
Let $\Theta$ be an $l$-dimensional face of an $n$-dimensional reflexive
polyhedron $\Delta \subset M_{\bf Q}$, $p_1, \ldots, p_k$ are vertices of
$\Theta$,
$\Delta^* \in N_{\bf Q}$ the dual reflexive polyhedron.

Define the dual to $\Theta$ $(n -l -1)$-dimensional face of $\Delta^*$ as
\[ \Theta^* = \{ y \in \Delta^* \mid \langle p_1, y \rangle  = \cdots =
\langle p_k, y \rangle = - 1 \}.  \]

Then one gets the one-to-one
correspondence $\Theta \leftrightarrow \Theta^*$  between faces of
the polyhedra $\Delta$ and $\Delta^*$ reversing the incidence relation of the
faces.
\label{dual.edge}
\end{prop}

 \begin{opr}
{\rm A complex normal irreducible $n$-dimensional projective
algebraic variety $W$ with only
Gorenstein canonical singularities we call a {\em Calabi-Yau variety}
if $W$ has trivial canonical
bundle and
\[ H^i (W, {\cal O}_W) = 0 \]
for $0 < i < n$.}
\label{calabi-yau}
\end{opr}

   The next theorem describes the relationship between reflexive
polyhedra and Calabi-Yau hypersurfaces.

\begin{theo}
Let $\Delta$ be an $n$-dimensional
integral polyhedron in $M_{\bf Q}$, ${\bf P}_{\Delta}$ the corresponding
$n$-dimensional projective toric variety, ${\cal F}(\Delta)$ the family
of projective $\Delta$-regular hypersurfaces $\overline{Z}_f$ in
${\bf P}_{\Delta}$.
Then the following conditions are equivalent

{\rm (i)} the family ${\cal F}(\Delta)$ of $\Delta$-hypersurfaces in
${\bf P}_{\Delta}$ consists of Calabi-Yau varieties with
canonical singularities $($see {\rm \ref{calabi-yau}}$);$

{\rm (ii)} the ample invertible sheaf ${\cal O}_{{\Delta}}(1)$
on the toric variety ${\bf P}_{\Delta}$ is anticanonical, i.e.,
${\bf P}_{\Delta}$ is a toric Fano variety with Gorenstein
singularities;

{\rm (iii)} $\Delta$ contains only one integral point $m_0$ in its
interior, and $(\Delta-m_0, M)$ is a reflexive pair.
\label{equiv}
\end{theo}

\proof  Since $\overline{Z}_f$ is an ample Cartier divisor on
${\bf P}_{\Delta}$, (i)$\Rightarrow$(ii) follows from the adjunction
formula.

The equivalence  (ii)$\Leftrightarrow$(iii)
follows from \ref{smooth.fano}.

Assume that (ii) and (iii) are satisfied. Let us prove (i).
By the adjunction formula, it follows from (ii) that every hypersurface
$\overline{Z}_f$ has
trivial canonical divisor. By the vanishing theorem for
arbitrary ample divisors on toric varieties \cite{dan1}, one  gets
\[ H^i(\overline{Z}_f, {\cal O}_{\overline{Z}_f}) = 0\]
for $ 0 < i < n-1$.
By \ref{anal2}, singularities of  $\overline{Z}_f$ are
analytically isomorphic to toric singularities of ${\bf P}_{\Delta}$.
Since all  singularities of ${\bf P}_{\Delta}$ are Gorenstein, by
\ref{gor.can}, they are also  canonical. So every $\Delta$-regular
hypersurface satisfies \ref{calabi-yau}.
\bigskip

The next statement follows  from definitions of the polyhedra
$\Delta(\Sigma, h_K)$ and $\Delta^*(\Sigma, h_K)$.

\begin{prop}
Let $(\Delta,M)$ be a reflexive pair, $(\Delta^*,N)$ the dual
reflexive pair, $\Sigma$ the rational polyhedral fan
defining the corresponding Gorenstein toric Fano variety ${\bf P}_{\Delta}$.
Then
\[ \Delta(\Sigma, h_K) = \Delta, \]
\[ \Delta^*(\Sigma, h_K) = \Delta^*. \]
In particular, if $\Sigma^*$ is a rational polyhedral fan defining
the Gorenstein toric Fano variety ${\bf P}_{\Delta^*}$, then
\[ \Delta(\Sigma, h_K) = \Delta^*(\Sigma^*, h_K) \]
and
\[ \Delta(\Sigma^*, h_K) = \Delta^*(\Sigma, h_K). \]
\end{prop}

Thus, in order to construct a rational polyhedral fan $\Sigma(\Delta)$
corresponding to a reflexive polyhedron $\Delta$, one can use the following
another way: we take the dual reflexive polyhedron $\Delta^*$ and apply
the statement

\begin{coro}
Let $\Delta \subset M_{\bf Q}$ be a reflexive polyhedron, $\Delta^*$
the dual reflexive polyhedron in $N_{\bf Q}$. For every
$l$-dimensional face $\Theta$ of $\Delta^*$ define the $(l+1)$-dimensional
convex cone $\sigma[\Theta]$ supporting the face  $\Theta$ as follows
\[ \sigma[\Theta] = \{ \lambda x \in M_{\bf Q} \mid
x \in \Theta,\;\; \lambda \in {\bf Q}_{\geq 0} \}. \]
Then the set $\Sigma[\Delta^*]$ of all cones $\sigma[\Theta]$ where
$\Theta$ runs over all faces of $\Delta^*$ is the complete fan defining
the toric Fano variety ${\bf P}_{\Delta}$ associated with $\Delta$.
Moreover, every
$(l+1)$-dimensional cone $\sigma[\Theta]$ coincides with
$\sigma(\Theta^*)$
{\rm ($see$ \ref{def.fan})}, where $\Theta^* \subset \Delta$ is the dual
to $\Theta$ $(n-l-1)$-dimensional face of $\Delta$
{\rm ($see$ \ref{dual.edge})}.
\end{coro}
\bigskip

There is the following finiteness theorem for reflexive polyhedra.

\begin{theo} There exist up to an unimodular transformation of
the lattice $M$ only finitely many reflexive pairs
$(\Delta, M)$ of fixed dimension $n$.
\end{theo}

This statement follows from the finiteness theorem in \cite{bat01}, or
from results in \cite{boris,hensley}.

\subsection{Singularities and morphisms of Calabi-Yau hypersurfaces}

\hspace*{\parindent}

Let $\Delta$ be a reflexive polyhedron, $\Delta^*$ the dual reflexive
polyhedron. Take a maximal projective triangulation ${\cal T}$ of $\Delta$.
It follows from the proof of \ref{crep.fano} that ${\cal T}$ defines a
$MPCP$-desingularization
\[ \varphi_{\cal T}\; :\; {\hat{\bf P}}_{\Delta} \rightarrow
{\bf P}_{\Delta} \]
of the Gorenstein toric variety ${\bf P}_{\Delta}$. Let $\overline{Z}_f$
be a $\Delta$-regular Calabi-Yau hypersurface in ${\bf P}_{\Delta}$.
Put
\[\hat{Z}_f = \varphi^{-1}_{\cal T}(\overline{Z}_f).\]
By \ref{max.crep},
\[ \varphi_{\cal T}\; : \;
 \hat{Z}_f \rightarrow \overline{Z}_f \]
 is a $MPCP$-desingularization of
$\overline{Z}_f$.

\begin{opr}
{\rm We will call
\[ \varphi_{\cal T}\; : \;  \hat{Z}_f \rightarrow \overline{Z}_f \]
the {\em toroidal  MPCP-desingularization of $\overline{Z}_f$ corresponding to
a maximal projective triangulation ${\cal T}$ of $\Delta^*$}.}
\end{opr}

Using \ref{crep.fano} and  \ref{codim2,3}(ii), one gets the following.

\begin{theo}
There exist at least one toroidal MPCP-desingularization $\hat{Z}_f$ of any
$\Delta$-regular Calabi-Yau hypersurface in ${\bf P}_{\Delta}$. Such a
MPCP-desingularization $\hat{Z}_f$ corresponds to any  maximal
projective triangulation ${\cal T}$ of the dual polyhedron $\Delta^*$.
The codimension of singularities of $\hat{Z}_f$ is always at least $4$.
\label{smooth.c}
\end{theo}

\begin{coro}
A toroidal MPCP-desingularization of a projective Calabi-Yau hypersurace
$\overline{Z}_f$ associated with a  reflexive polyhedron
$\Delta$ of dimension $n \leq 4$ is always a  smooth Calabi-Yau manifold.
\label{smooth.c1}
\end{coro}

Let ${\cal T}$ be a maximal projective triangulation of $\Delta^*$. For
any  $l$-dimensional face $\Theta$ of $\Delta$, the restriction of
${\cal T}$ on the dual $(n - l -1)$-dimensional face $\Theta^* \subset
\Delta^*$ is a maximal projective triangulation ${\cal T}\mid_{\Theta^*}$ of
$\Theta^*$. By \ref{anal}, the analytical decription  of singularities along
a stratum $Z_{f,\Theta}$  as well as of their  $MPCP$-desingularizations
reduces to
the combinatorial description of a $MPCP$-desingularization of
the toric singularity at the unique
closed ${\bf T}_{\sigma}$-invariant point $p_{\sigma}$  on
the $(n-l)$-dimensional affine toric variety ${\bf A}_{\sigma, N(\sigma)}$,
 where $\sigma =  \sigma[\Theta^*]$. So we introduce the following definition.

\begin{opr}
{\rm  We call the face $\Theta^*$ of the polyhedron $\Delta$
the {\em diagram} of the toric singularity
at $p_{\sigma} \in  {\bf A}_{\sigma, N(\sigma)}$. A maximal projective
triangulation ${\cal T}\mid_{\Theta^*}$ of $\Theta^*$ induced by a
maximal projective triangulation ${\cal T}$ of $\Delta^*$
we call {\em a triangulated diagram}  of the toric singularity at
$p_{\sigma} \in  {\bf A}_{\sigma, N(\sigma)}$.   }
\label{diagram}
\end{opr}

We have seen already in \ref{anal2} that for any face $\Theta \subset \Delta$
and
any $\Delta$-regular Laurent polynomial $f \in L(\Delta)$
local neighbourhoods of points belonging to the same stratum $Z_{f, \Theta}$
are analitically  isomorphic. Thus, if a stratum $Z_{f, \Theta}$ consists of
singular points of $\overline{Z}_f$, then all these singularities are
analitically isomorphic. Our purpose now is to describe singularities along
$Z_{f, \Theta}$ and their $MPCP$-desingularizations in terms of triangulated
diagrams.
\medskip

Let ${\cal T}\mid_{\Theta^*}$ be a triangulated diagram. Then we obtain
a subdivision $\Sigma({\cal T},\Theta)$ of the cone $\sigma = \sigma[\Theta^*]$
into the union of subcones supporting elementary simplices of the
triangulation  ${\cal T}\mid_{\Theta^*}$. By \ref{biration}, one has
the corresponding  projective birational toric morphism
\[ \varphi_{{\cal T}, \Theta^*}\; : \;{\bf P}_{\Sigma({\cal T},\Theta)}
\rightarrow {\bf A}_{\sigma, N(\sigma)}. \]

\begin{theo}
For any $l$-dimensional face $\Theta \subset \Delta$ and any
closed point $p \in Z_{f, \Theta}$, the fiber $\varphi_{\cal T}^{-1}(p)$ of
a MPCP-desingularization $\varphi_{\cal T}$ is isomorphic to the fiber
$\varphi_{{\cal T}, \Theta^*}(p_{\sigma})$ of the projective toric morphism
$\varphi_{{\cal T}, \Theta^*}$.

The number of irreducible $(n-l-1)$-dimensional components of
$\varphi_{\cal T}^{-1}(p)$ equals $l^*(\Theta^*)$, i.e., the number of
integral points in the interior of $\Theta^*$. Moreover,
the Euler characteristic
of $\varphi_{\cal T}^{-1}(p)$ equals the number of elementary simplices
in the triangulated diagram ${\cal T}\mid_{\Theta^*}$.
\label{topol.des}
\end{theo}

\proof  Since $ \varphi_{\cal T}\; :\; {\hat{\bf P}}_{\Delta} \rightarrow
{\bf P}_{\Delta}$ is a birational toric morphism, the ${\bf T}$-action
induces isomorphisms of fibers of $\varphi_{\cal T}$ over closed points of a
${\bf T}$-stratum  ${\bf T}_{\Theta}$.
Thus, we obtain isomorphisms between fibers
of $\varphi_{\cal T}$ over closed points of  $Z_{f, \Theta} \subset
{\bf T}_{\Theta}$. By \ref{orbits.fan}(i),
${\bf T}_{\Theta}$ is contained in the
${\bf T}$-invariant open subset ${\bf A}_{\sigma,N} \cong {\bf T}_{\Theta}
\times {\bf A}_{\sigma, N(\sigma)}$. Thus, we have a commutative
diagram
\[
\begin{tabular}{ccc}
$ \varphi_{\cal T}^{-1}({\bf A}_{\sigma,N})$ &  $\rightarrow$ & $
{\bf A}_{\sigma,N}$ \\
$\downarrow$ &   & $\downarrow$ \\
$ {\bf P}_{\Sigma({\cal T},\Theta)}
$ & $\rightarrow$ & ${\bf A}_{\sigma, N(\sigma)}$
\end{tabular} \]
whose vertical maps are divisions by the action of the torus
${\bf T}_{\Theta}$.
So the fiber $\varphi_{\cal T}^{-1}(p)$ over any closed point $p \in
{\bf T}_{\Theta^*}$  is isomorphic to the fiber
$\varphi_{{\cal T}, \Theta^*}^{-1}(p_{\sigma})$ of the projective
toric morphism $\varphi_{{\cal T}, \Theta^*}$. Therefore,
irreducible divisors of
${\bf P}_{\Sigma({\cal T},\Theta)}$ over $p_{\sigma}$ one-to-one correspond
to integral points in the interior of $\Theta^*$.

On the other hand, the toric morphism $\varphi_{{\cal T}, \Theta^*}$ has
an action the $(n-l)$-dimensional torus ${\bf T}_{\Theta}' =
{\rm Ker}\, \lbrack {\bf T} \rightarrow {\bf T}_{\Theta} \rbrack$. Since the
closed point $p_{\sigma} \in {\bf A}_{\sigma, N(\sigma)}$ is
${\bf T}_{\Theta}'$-invariant, one has a decomposition of the fiber
$\varphi_{{\cal T}, \Theta^*}^{-1}(p_{\sigma})$ into the union of
${\bf T}_{\Theta}'$-orbits. Thus the Euler characteristic of
$\varphi_{{\cal T}, \Theta^*}^{-1}(p_{\sigma})$ is the number of
zero-dimensional ${\bf T}_{\Theta}'$-orbits. The latter equals the
number of $(n-l)$-dimensional cones of $\Sigma({\cal T},\Theta)$, i.e.,
the number of elementary simplices in the triangulated diagram
${\cal T}\mid_{\Theta^*}$.
\medskip

\begin{exam}
{\rm Let $\Theta$ be an $(n-2)$-dimensional face of a reflexive polyhedron
$\Delta$, and let $\Theta^*$ be the dual to $\Theta$ $1$-dimensional face
of the dual
polyhedron $\Delta^*$. There exists a unique maximal projective
triangulation of $\Theta^*$ consisting of $d(\Theta^*)$ elementary,
in fact,  regular segments.

In this case, small analytical
neighbourhoods of points on $Z_{f, \Theta}$ are analytically isomorphic to
product of $(n-3)$-dimensional open ball and a small analytical
neighbourhood of $2$-dimensional double point singularity of type
$A_{d(\Theta^*)-1}$.

The fiber of $\varphi_{\cal T}$ over any point of $\overline{Z}_{f,\Theta}$
is the Hirzebruch-Jung tree of $l^*(\Theta^*) = d(\Theta^*) -1$
smooth rational curves having
an action of ${\bf C}^*$. }
\label{duval}
\end{exam}

\begin{opr}
{\rm Let $(\Delta_1, M_1)$ and $(\Delta_2, M_2)$ be two reflexive
pairs of equal dimension. {\em A finite morphism} of reflexive pair
\[ \phi\; :\; (\Delta_1, M_1) \rightarrow (\Delta_2, M_2) \]
is  a homomorphism of lattices $\phi\,:\,  M_1 \rightarrow M_2$ such that
$\phi (\Delta_1) = \Delta_2$.
}
\label{morph}
\end{opr}

By \ref{final}, we obtain:

\begin{prop}
Let
\[ \phi\; :\; (\Delta_1, M_1) \rightarrow (\Delta_2, M_2) \]
be a finite morphism of reflexive pairs. Then the dual finite morphism
\[ \phi^* \; : \; (\Delta_2^*, N_2) \rightarrow (\Delta_1^*, N_1). \]
induces a finite  surjective morphism
\[ \tilde {\phi}^* \; : \; {\bf P}_{\Delta_2,M_2}
\rightarrow {\bf P}_{\Delta_1,M_1}.  \]
Moreover, the restriction
\[ \tilde {\phi}^* \; : \; {\bf P}_{\Delta_2,M_2}^{[1]}
\rightarrow {\bf P}_{\Delta_1,M_1}^{[1]}  \]
is an \^etale morphism of degree
\[ d_{M_2}(\Delta_2) / d_{M_1}(\Delta_1) =
d_{N_1}(\Delta_1^*) / d_{N_2}(\Delta_2^*).\]
\label{etale}
\end{prop}

\proof  It remains to show that
$\tilde {\phi}^* \; : \; {\bf P}_{\Delta_2,M_2}^{[1]}
\rightarrow {\bf P}_{\Delta_1,M_1}^{[1]}$ is  \^etale.
Take a $\Delta_1$-regular Calabi-Yau hypersurface $\overline{Z}_{f}
\subset {\bf P}_{\Delta_1,M_1}^{[1]}$. By \ref{subdiv.hyp}, the
hypersurface $\overline{Z}_{\tilde{\phi}^*f} =
\tilde{\phi}^{-1}(\overline{Z}_f)$ is $\Delta_2$-regular.
By \ref{sing.hyp} and \ref{equiv}, two quasi-projective  varieties
$Z^{[1]}_f$ and
${Z}_{\tilde{\phi}^*f}^{[1]}$ are smooth and have trivial canonical
class. Therefore, any finite morphism of these varieties must be
\^etale.

\subsection{The Hodge  number $h^{n-2,1}(\hat{Z}_f)$}

\hspace*{\parindent}

Let us  apply the result of Danilov and Khovanski\^i (see \ref{hd.dh})
to  the case of reflexive polyhedra $\Delta$ of dimension $\geq 4$.
Using the properties $l^*(2\Delta) = l(\Delta)$, $l^*(\Delta) =1$, we
can calculate  the Hodge-Deligne number
$h^{n-2,1}(Z_f)$ of an affine Calabi-Yau hypersurface
$Z_f$ as follows.

\begin{theo}
Let $\Delta$ be a reflexive $n$-dimensional polyhedron $(n\geq 4)$, then
the  Hodge-Deligne number $h^{n-2,1}$ of the cohomology group $H^{n-1}(Z_f)$
of any $(n-1)$-dimensional affine $\Delta$-regular Calabi-Yau
hypersurface $Z_f$ equals
\[ h^{n-2,1} (Z_f) = l(\Delta)  - n - 1 -  \sum_{{\rm codim}\, \Theta =1 }
l^*(\Theta). \]
\label{n.def}
\end{theo}
\bigskip

In fact, we can calculate the Hodge-Deligne space $H^{n-2,1}(Z_f)$ itself
(not only the dimension $h^{n-2,1}$).

\begin{theo}
{\rm \cite{batyrev.var}}
Let $L^*_1(\Delta)$ be the subspace in $L(\Delta)$ generated by all
monomials $X^m$
such  that $m$ is an interior integral point on a face
$\Theta \subset \Delta$ of codimension $1$. Then the ${\bf C}$-space
$H^{n-2,1}(Z_f)$ is canonically isomorphic to the quotient of $L(\Delta)$ by
\[ L^*_1(\Delta) + {\bf C}\langle f, f_1, \ldots, f_n \rangle,  \]
where
\[ f_i (X) = X_i \frac{\partial }{\partial X_i} f(X),
\;\; (1 \leq i \leq n). \]
\end{theo}

We want now to calculate for $n \geq 4$ the Hodge number $h^{n-2,1}$ of a
{\em MPCP}-de\-sin\-gu\-la\-ri\-za\-tion ${\hat Z}_f$ of the toroidal
compactification
$\overline{Z}_f$ of $Z_f$.

For this purpose, it is  convenient to use the notion of the
$(p,q)$-{\em Euler characteristic} $e_c^{p,q}$ introduced in \cite{dan.hov}.

\begin{opr}
{\rm For any complex algebraic variety $V$, $e_c^{p,q}(V)$ is defined as the
alternated sum of Hodge-Deligne numbers
\[ \sum_{i \geq 0} (-1)^i h^{p,q}(H^i_c(V)). \] }
\end{opr}

\begin{prop}({\rm see \cite{dan.hov}})
Let  $V = V' \times V''$ be a product of two
algebraic varieties. Then one has
\[ e_c^{p,q}(V) = \sum_{(p'+ p'',q' +q'') = (p,q)}
e_c^{p',q'}(V') \cdot e_c^{p'',q''}(V''). \]
\label{h.product}
\end{prop}

\begin{opr}
{\rm A stratification of a compact algebraic variety $V$ is a representation
of $V$ as a disjoint union of finitely many locally closed smooth
subvarieties $\{ V_j \}_{j \in J}$ (which are called {\em strata})
such that for any $j \in J$ the closure of $V_j$ in $V$ is a union of
the strata.}
\end{opr}

The following property follows immediately from long cohomology sequences
(see \cite{dan.hov}).

\begin{prop}
Let $\{ V_j \}_{j \in J}$ be a stratification of $V$. Then
\[ e_c^{p,q}(V) = \sum_{j \in J} e_c^{p,q}(V_j). \]
\label{addit}
\end{prop}

Returning to our Calabi-Yau variety $\hat{Z}_f$, we see that $\hat{Z}_f$
is always quasi-smooth. Therefore the cohomology groups $H_c^i(\hat{Z}_f)
\cong H^i(\hat{Z}_f)$ have the pure Hodge structure of weight $i$. So,
by \ref{addit},
it suffices to calculate the $(n-2,1)$-Euler characteristic
\[ e_c^{n-2,1}(\hat{Z}_f) = (-1)^{n-1} h^{n-2,1}(\hat{Z}_f). \]
\medskip

First, we define a convenient stratification of $\hat{Z}_f$.

Let $\varphi_{\cal T} \,: \, \hat{Z}_f \rightarrow
\overline{Z}_f$ be the corresponding birational morphism. Then
$\hat{Z}_f$ can be represented as a disjoint union
\[ \hat{Z}_f = \bigcup_{\Theta \subset \Delta}
\varphi^{-1}_{\cal T}(Z_{f, \Theta}). \]
On the other hand, all irreducible components of fibers of $\varphi_{\cal T}$
over closed points of  $Z_{f, \Theta}$ are toric varieties.
Therefore, we can define
a stratification of $\varphi^{-1}_{\cal T}(Z_{f, \Theta})$ by smooth
affine  algebraic varieties which are isomorphic to
products ${Z}_{f,\Theta} \times ({\bf C}^*)^k$ some nonnegative integer $k$.
As a result, we obtain a stratification of  $\hat{Z}_f$
by  smooth affine varieties which are isomorphic to
products ${Z}_{f,\Theta} \times ({\bf C}^*)^k$ for some face $\Theta \subset
\Delta$ and for some nonnegative integer $k$.

Second, we note that $(n-2,1)$-Euler characteristic of
${Z}_{f,\Theta} \times ({\bf C}^*)^k$ might be nonzero only in two cases:
$\Theta = \Delta$, or ${\rm dim}\, \Theta = n -2$ and $k =1$. The latter
follows from \ref{properties.hodge}, \ref{h.product}, and from the observation
that the $(p,q)$-Euler characteristic of an
algebraic torus $({\bf C}^*)^k$ is nonzero only if $p =q$ .

We already know from \ref{n.def} that
\[ e_c^{n-2,1}(Z_f) = (-1)^{n-1}(l(\Delta)  - n - 1 -
\sum_{{\rm codim}\, \Theta =1 } l^*(\Theta)). \]
On the other hand, the strata which are isomorphic to
${Z}_{f,\Theta} \times {\bf C}^*$ appear from the fibers of $\varphi_{\cal T}$
over $(n-3)$-dimensional singular affine locally closed subvarieites
${Z}_{f,\Theta} \subset \overline{Z}_f$ having codimension $2$ in
$\overline{Z}_f$.
By \ref{duval}, a $\varphi^{-1}(Z_{f, \Theta})$ consists of $l^*(\Theta^*) =
d(\Theta^*) -1$ irreducible components  isomorphic to ${\bf P}_1 \times
{Z}_{f,\Theta}$. As a result, for every $(n-2)$-dimensional face $\Theta
\subset \Delta$, one obtains  $l^*(\Theta^*)$  strata isomorphic to
${Z}_{f,\Theta} \times {\bf C}^*$.
On the other hand, one has

\[ e_c^{n-2,1}({Z}_{f,\Theta} \times {\bf C}^*) =
e_c^{n-3,0}({Z}_{f,\Theta}) \cdot e_c^{1,1}({\bf C}^*). \]

It is clear that $e_c^{1,1}({\bf C}^*) =1$.
By results of Danilov and Khovanski\^i (see \cite{dan.hov}), one has
\[ e_c^{n-3,0}({Z}_{f,\Theta}) = (-1)^{n-3} l^*(\Theta). \]
Thus we come to the following result.

\begin{theo}
For $n \geq 4$, the Hodge number $h^{n-2,1}({\hat Z}_f)$ equals
\[ l(\Delta)  - n - 1 -
\sum_{{\rm codim}\, \Theta =1 }
l^*(\Theta) + \sum_{{\rm codim}\, \Theta =2 }
l^*(\Theta) \cdot l^*(\Theta^*), \]
where $\Theta$ denotes a face of a reflexive $n$-dimensional
polyhedron $\Delta$, and $\Theta^*$ denotes the corresponding dual
face of the dual reflexive polyhedron $\Delta^*$.
\label{n.def.com}
\end{theo}

\begin{coro}
Assume that ${\bf P}_{\Delta}$, or ${\bf P}_{\Delta^*}$ is a smooth
toric Fano variety of dimension $n \geq 4$ (see \ref{smooth.fano}(ii)).
Then the Hodge-Deligne number $h^{n-2,1}$ of an affine
$\Delta$-regular hypersurface $Z_f$ coincides with the Hodge number
of a MPCP-desingularization $\hat{Z}_f$ of its toroidal
compactification $\overline{Z}_f$.
\end{coro}

\proof If ${\bf P}_{\Delta}$, or ${\bf P}_{\Delta^*}$ is a smooth
toric Fano variety, then for any face $\Theta$ of codimension $2$,
one has $l^*(\Theta)=0$, or $l^*(\Theta^*)=0$.
\bigskip

\subsection{The Hodge number $h^{1,1}(\hat{Z}_f)$}

\hspace*{\parindent}

First we note that the group of principal ${\bf T}$-invariant
divisors is isomorphic to the lattice $M$.
Applying  \ref{crep.fano}, we obtain:

\begin{prop}
Any  MPCP-desingularization $\hat{\bf P}_{\Delta}$ of a toric Fano variety
${\bf P}_{\Delta}$ contains exactly
\[ l(\Delta^*) - 1 = {\rm card}\, \{ N \cap \partial \Delta^* \} \]
${\bf T}$-invariant divisors, i.e., the Picard number
$\rho(\hat{\bf P}_{\Delta}) =  h^{1,1}(\hat{\bf P}_{\Delta})$
equals
\[   l(\Delta^*) - n - 1. \]
\label{bound.points}
\end{prop}

\begin{theo}
Let $\hat{Z}_f$ be a MPCP-desingularization
of a projective $\Delta$-regular Calabi-Yau hypersurface
$\overline{Z}_f$, then for $n \geq 4$ the Hodge number
$h^{1,1}(\hat{Z}_f)$, or the Picard number of $\hat{Z}_f$, equals
\begin{equation}
l(\Delta^*) - n - 1 -  \sum_{{\rm codim}\, \Theta^* =1 } l^*(\Theta^*) +
 \sum_{{\rm codim}\, \Theta =2 } l^*(\Theta^*) \cdot l^*(\Theta),
\label{for.pic}
\end{equation}
where $\Theta^*$ denotes a face of the dual to $\Delta$
reflexive $n$-dimensional
polyhedron $\Delta^*$, and $\Theta$ denotes the corresponding dual to $\Theta$
face of the  reflexive polyhedron $\Delta$.
\label{n.pic}
\end{theo}

\proof
By \ref{crep.fano},
a $MPCP$-desingularization $\hat{Z}_f$ of a  $\Delta$-regular
Calabi-Yau hypersurface $\overline{Z}_f$ is  induced by
a $MPCP$-desingularization $\varphi_{\cal T}\,: \,
\hat{\bf P}_{\Delta} \rightarrow {\bf P}_{\Delta}$ of the ambient
toric Fano variety ${\bf P}_{\Delta}$.
Since $\hat{Z}_f$ has only terminal
${\bf Q}$-factorial singularities, i.e., $\hat{Z}_f$ is
quasi-smooth \cite{dan1},  the Hodge structure in cohomology groups
of $\hat{Z}_f$ is pure and satisfies the Poincare duality
(see \ref{duality}).
Therefore, the number $h^{1,1}$ equals the Hodge number $h^{n-2,n-2}$ in the
cohomology group $H^{n-2}_c (\hat{Z}_f)$ with compact
supports. By the Lefschetz-type theorem (see \ref{Lef}), if $n \geq 4$, then
for $i = n-3, n-4$ the Gysin homomorphisms
\[ H^i_c(Z_f) \rightarrow H^{i+2}_c ({\bf T}) \]
are isomorphisms  of Hodge structures with the shifting
the Hodge type by $(1,1)$. On the other hand,  $H^{2n-1}_c ({\bf T})$
is an $n$-dimensional space having the Hodge
type $(n-1,n-1)$, and the space $H^{2n-2}_c ({\bf T})$ has the Hodge
type $(n-2,n-2)$. So
$H^{2n-3}_c (Z_f)$ is an $n$-dimensional space having the Hodge type
$(n-2,n-2)$, and the space $H^{2n-4}_c (Z_f)$ has the Hodge type $(n-3,n-3)$.
The complementary set $Y = \hat{Z}_f \setminus Z_f$ is a
closed subvariety of $\hat{Z}_f$ of codimension $1$. Consider
the corresponding exact sequence of Hodge structures
\[ \cdots \rightarrow H_c^{2n-4}(Z_f) \stackrel{\beta_1}{\rightarrow}
H_c^{2n-4}(\hat{Z}_f) {\rightarrow}
H_c^{2n-4}(Y) {\rightarrow}  \\
H_c^{2n-3}(Z_f) \stackrel{\beta_2}{\rightarrow}
H_c^{2n-3}(\hat{Z}_f) {\rightarrow} \cdots \]
Comparing the Hodge types, we immediately get that $\beta_1$ and $\beta_2$ are
zero mappings.
Since the space $H_c^{2n-4}(Y)$ does not have subspaces of the Hodge type
$(n-1, n-3)$, the Hodge humber $h^{n-2,n-2}$ of
$H_c^{2n-4}(\hat{Z}_f)$ equals
${\rm dim}\, H_c^{2n-4} (\hat{Z}_f)$.
Thus we get the short exact sequence of cohomology
groups of the Hodge type $(n-2, n-2)$
\[ 0 \rightarrow  H_c^{2n-4}(\hat{Z}_f) {\rightarrow}
H_c^{2n-4}(Y) {\rightarrow}
H_c^{2n-3}(Z_f) {\rightarrow} 0.\]

It is easy to see that the dimension of $H_c^{2n-4}(Y)$ equals the number
of the irreducible components of the $(n-2)$-dimensional complex subvariety
$Y$. On the other hand, $Y$ is
the intersection of $\hat{Z}_f \subset \hat{\bf P}_{\Delta}$ with the
union of all irreducible
${\bf T}$-invariant  divisors
on the corresponding maximal partial
crepant desingularization $\hat{\bf P}_{\Delta}$
of the toric Fano variety ${\bf P}_{\Delta}$.
We have seen in
\ref{bound.points} that these  toric
divisors on $\hat{\bf P}_{\Delta}$ are in the one-to-one correspondence to the
integral points $\rho \in N \cap \partial \Delta^*$,i.e., we have
exactly $l(\Delta^*) -1$ irreducible
toric divisors on $\hat{\bf P}_{\Delta}$. Every such a divisor $D_{\rho}$
is the closure of a $(n-1)$-dimensional torus ${\bf T}_{\rho}$
whose lattice of characters consists of elements of $M$ which are orthogonal
to $\rho \in N$. It is important to note that $D_{\rho} \cap \hat{Z}_f$ is
the closure in $\hat{\bf P}_{\Delta}$ of the affine hypersurface
$\varphi_{\cal T}^{-1}(Z_{f,\Theta} \cap {\bf T}_{\rho} \subset
{\bf T}_{\rho}$ where $\varphi_{\cal T}({\bf T}_{\rho}) = {\bf T}_{\Theta}$.

Note that since
$\overline{Z}_f$
does not intersect any $0$-dimensional ${\bf T}$-orbit (\ref{zero.orb}),
$\hat{Z}_f$ does not intersect  exeptional divisors on
$\hat{\bf P}_{\Delta}$ lying over these points. The exceptional divisors
lying over ${\bf T}$-invariant points of ${\bf P}_{\Delta}$ correspond
to  integral points $\rho$ of  $N$ in the interiors of $(n-1)$-dimensional
faces $\Theta^*$ of $\Delta^*$. So we must consider only
\[ l(\Delta^*) -1 - \sum_{{\rm codim}\, \Theta^* =1 }
l^*(\Theta^*) \]
integral points of $N \cap \partial \Delta$ which are contained in faces
of codimension $\geq 2$.

If $D_{\rho}$ is an invariant  toric divisor on $\hat{\bf P}_{\Delta}$
corresponding to an integral point $\rho$ belonging to the interior of a
$(n-2)$-dimensional face $\Theta^*$ of $\Delta^*$, then $D_{\rho} \cap
\hat{Z}_f$ consists of $d(\Theta)$ irreducible components whose
$\varphi_{\cal T}$-images are $d(\Theta)$ distict points of the
zero-dimensional  stratum $Z_{f, \Theta} \subset {\bf T}_{\Theta}$.

If $D_{\rho}$ is an invariant  toric divisor on $\hat{\bf P}_{\Delta}$
corresponding to an integral point belonging to the interior of a
face $\Theta^* \subset \Delta^*$ of codimension $\geq 3$, then $D_{\rho} \cap
\hat{Z}_f$ is  irreducible because ${\bf T}_{\rho} \cap \hat{Z}_f$ is
an irreducible affine hypersurface in ${\bf T}_{\rho}$ isomorphic
to $Z_{f, \Theta} \times ({\bf C}^*)^{n-1 - {\rm dim}\,\Theta}$.

Consequently, the number of irreducible components of $Y$ is
\[ \sum_{{\rm codim}\, \Theta^* =2 } d(\Theta)\cdot l^*(\Theta^*) + \]
\[ + {\rm number\; of\; integral\; points\; on\; faces}\;
\Theta^* \subset \Delta^*, \; {\rm codim}\, \Theta^* \geq 3. \]

Since $d(\Theta) = l^*(\Theta) +1$ for any  $1$-dimensional face
$\Theta \Delta$, we can rewrite the number of the irreducible
components of $Y$  as follows
\[ {\rm dim}\, H_c^{2n-4}(Y) = l(\Delta^*) - 1 -
\sum_{{\rm codim}\, \Theta^* =1 }
l^*(\Theta^*) +  \sum_{{\rm codim}\, \Theta =2 }
l^*(\Theta^*) \cdot l^*(\Theta). \]
Since  ${\rm dim}\, H_c^{2n-3}(Z_f) = n$, we obtain $(\ref{for.pic})$.
\bigskip

Applying  \ref{n.def.com} and \ref{n.pic}, we conclude.

\begin{theo}
For any reflexive polyhedron $\Delta$ of dimension $n \geq 4$,
the Hodge number $h^{n-1,1}(\hat{Z}_f)$ of a MPCP-desingularization
of a $\Delta$-regular Calabi-Yau hypersurface $\overline{Z}_f \subset
{\bf P}_{\Delta}$ equals the
Picard  number $h^{1,1}(\hat{Z}_g)$  of  a MPCP-desingularization
of a $\Delta^*$-regular projective Calabi-Yau
hypersurface $\overline{Z}_g \subset {\bf P}_{\Delta^*}$ corresponding to
the dual reflexive polyhedron $\Delta^*$.
\label{num.mir}
\end{theo}

\subsection{Calabi-Yau $3$-folds}

\hspace*{\parindent}

We begin with the remark that  only if ${\rm dim}\, \Delta =4$
both statements \ref{smooth.c1} and  \ref{num.mir} hold.
In this case,  we deal with  $3$-dimensional Calabi-Yau hypersurfaces
$\overline{Z}_f \subset {\bf P}_{\Delta}$ which admit   {\em smooth}
$MPCP$-desingularizations $\hat{Z}_f$.

 Calabi-Yau $3$-folds $\hat{Z}_f$  are of primary interest in
theoretical  physics. The number
\[ \frac{1}{2} e(\hat{Z}_f) =
( h^{1,1}(\hat{Z}_f) - h^{2,1}(\hat{Z}_f)) \]
is called {\em the number
of generations} in superstring theory \cite{witt}.
So it is important to have
a simple formula for the Euler characteristic $e(\hat{Z}_f)$.

We have already calculated the Hodge numbers  $h^{1,1}(\hat{Z}_f)$ and
$h^{n-2,1}({Z}_f)$. So we obtain

\begin{coro}
For any Calabi-Yau $3$-fold $\hat{Z}_f$ defined by a $\Delta$-regular
Laurent polynomial $f$ whose Newton polyhedron is a reflexive
$4$-dimensional polyhedron $\Delta$, one has the following formulas
for the Hodge numbers
\[ h^{1,1}(\hat{Z}_f) =
l(\Delta^*) - 5 -  \sum_{{\rm codim}\, \Theta^* =1 } l^*(\Theta^*) +
 \sum_{{\rm codim}\, \Theta =2 } l^*(\Theta^*) \cdot l^*(\Theta), \]
\[ _+ h^{2,1} (\hat{Z}_f) = l(\Delta)  - 5 -  \sum_{{\rm codim}\, \Theta =1 }
l^*(\Theta) + \sum_{{\rm codim}\, \Theta =2 } l^*(\Theta) \cdot
l^*(\Theta^*). \]
\label{euler.iso}
\end{coro}

This implies also.

\begin{coro}
\[ e(\hat{Z}_f) = (l(\Delta)  - l(\Delta^*)) -
\; \left( \sum_{ \begin{array}{c} {\scriptstyle {\rm codim}\, \Theta =1} \\
{\scriptstyle \Theta \subset \Delta} \end{array}} l^*(\Theta)  -
\sum_{\begin{array}{c} {\scriptstyle {\rm codim}\, \Xi =1} \\
{\scriptstyle \Xi \subset \Delta^* } \end{array}} l^*(\Xi)\;\; \right) + \]
\[ + \left( \sum_{\begin{array}{c} {\scriptstyle {\rm codim}\, \Theta =2}  \\
{\scriptstyle \Theta \subset \Delta} \end{array}} l^*(\Theta) \cdot
l^*(\Theta^*) -
\sum_{\begin{array}{c}{\scriptstyle {\rm codim}\, \Xi =2} \\
{\scriptstyle \Xi \subset \Delta^*} \end{array}}
l^*(\Xi) \cdot l^*(\Xi^*)\; \right).  \]
\end{coro}
\medskip

Now we prove a more simple another
formula for the Euler characteristic of
Calabi-Yau $3$-folds.

\begin{theo}
Let $\hat{Z}_f$ be a MPCP-desingularization of a $3$-dimensional
$\Delta$-regular Calabi-Yau hypersurface associated with a $4$-dimensional
reflexive polyhedron $\Delta$. Then
\[ e(\hat{Z}_f) =
\sum_{\begin{array}{c}{\scriptstyle {\rm dim}\, \Theta =1} \\
{\scriptstyle \Theta \subset \Delta} \end{array}} d(\Theta)d(\Theta^*)  -
\sum_{\begin{array}{c}{\scriptstyle {\rm dim}\, \Theta =2} \\
{\scriptstyle \Theta \subset \Delta} \end{array}} d(\Theta)d(\Theta^*). \]
\label{new.euler}
\end{theo}

In our proof of theorem \ref{new.euler} we will use
one  general property of smooth quasi-projective open
subsets ${Z}_f^{[1]}$ in   $\overline{Z}_f$ consisting of the
union of the affine part ${Z}_f$ and all affine strata $Z_{f,\Theta}$, where
$\Theta$ runs over all faces of $\Delta$ of codimension 1, i.e.,
\[{Z}^{[1]}_f = {\bf P}_{\Delta}^{[1]} \cap {\overline{Z}}_f . \]

\begin{theo}
For arbitrary $n$-dimensional reflexive polyhedron $\Delta$ and
$\Delta$-regular Laurent polynomial $f \in L(\Delta)$, the
Euler characteristic   of the smooth quasi-projective   Calabi-Yau
variety ${Z}_f^{[1]}$
is always zero.
\label{euler.zero}
\end{theo}

\proof  Since  ${Z}_f^{[1]}$ is smooth,  the Euler characteristic
of the usual cohomology groups $H^* ({Z}_f^{[1]})$  is zero if and only if
the Euler characteristic  of the cohomology groups with compact
supports $H_c^* ({Z}_f^{[1]})$ is zero (see \ref{duality}).
It follows from the long exact sequence
of cohomology groups with compact suport that
\[ e(H_c^* ({Z}_f^{[1]})) = e(H_c^* ({Z}_f)) +
\sum_{{\rm codim \, \Theta} = 1}
e(H_c^* ({Z}_{f,\Theta})). \]
By \ref{euler},
\[ e(H_c^* ({ Z}_f)) = (-1)^{n-1} d(\Delta), \;\;
e(H_c^* ({ Z}_{f,\Theta})) = (-1)^{n-2} d(\Theta). \]
Thus,  it is sufficient to prove
\[ d(\Delta) = \sum_{{\rm codim}\,  \Theta =1} d(\Theta). \]
The latter follows immediately from the representation of the
$n$-dimensional polyhedron $\Delta$ as a union of $n$-dimensional
pyramids with vertex $0$ over all
$(n-1)$-dimensional faces $\Theta \subset \Delta$.
\bigskip

{\bf Proof of Theorem \ref{new.euler}}.
Let $\varphi_{\cal T}\;:\; \hat{Z}_f \rightarrow \overline{Z}_f$ be a
$MPCP$-desingularization. For any face $\Theta \subset \Delta$, let
$F_{\Theta}$ denotes $\varphi^{-1}_{\cal T}(Z_{f,\Theta})$.
We know that $\varphi$ is an isomorphism over
$Z^{[1]}_f$ which is the union of
the strata $Z_{f, \Theta}$ (${\rm dim}\, \Theta = 3,4$).
By  \ref{euler.zero}, $e(Z^{[1]}_f) =0$. Using additivity property of
the Euler characteristic, we obtain
\[ e(\hat{Z}_f) =
\sum_{\begin{array}{c}{\scriptstyle {\rm dim}\, \Theta =1} \\
{\scriptstyle \Theta \subset \Delta} \end{array}} e(F_{\Theta}) +
\sum_{\begin{array}{c}{\scriptstyle {\rm dim}\, \Theta =2} \\
{\scriptstyle \Theta \subset \Delta} \end{array}} e(F_{\Theta}). \]
Now the statement follows from the following lemma.

\begin{lem}
Let $\varphi_{\cal T}$ be a MPCP-desingularization as above. Then the
Euler characteristic $e(F_{\Theta})$
equals $d(\Theta)d(\Theta^*)$ if ${\rm dim}\, \Theta =1$, and
$-d(\Theta)d(\Theta^*)$ if ${\rm dim}\, \Theta =2$.
\end{lem}

\proof  If ${\rm dim}\, \Theta =1$, then $Z_{f,\Theta}$ consists of
$d(\Theta)$ distinct points. The fiber $\varphi^{-1}_{\cal T}(p)$
over any such a
point $p \in Z_{f,\Theta}$ is a union of smooth toric surfaces.
By \ref{topol.des},
the Euler characteristic of $\varphi^{-1}_{\cal T}(p)$ equals the number of
elementary simplices in the corresponding maximal projective
triangulation of $\Theta^*$. Since ${\rm dim}\, \Theta^* =2$, the
number of elementary simplices equals $d(\Theta^*)$ (see \ref{el.reg}). Thus,
$e(F_{\Theta}) = e(Z_{f,\Theta})e(\varphi^{-1}_{\cal T}(p)) =
d(\Theta)d(\Theta^*)$.

If ${\rm dim}\, \Theta =2$, then $Z_{f,\Theta}$ is a smooth affine algebraic
curve. By \ref{euler}, the Euler characteristic of $Z_{f,\Theta}$ equals
$-d(\Theta)$. By \ref{duval}, the fiber $\varphi^{-1}_{\cal T}(p)$ over
any point $p \in
Z_{f,\Theta}$ is the Hirzebruch-Jung tree  of $d(\Theta^*)-1$ smooth rational
curves, i.e., the Euler characteristic of $\varphi^{-1}_{\cal T}(p)$ again
equals $d(\Theta^*)$. Thus,
$e(F_{\Theta}) = -d(Z_{f,\Theta})e(\varphi^{-1}_{\cal T}(p)) =
-d(\Theta)d(\Theta^*)$.
\medskip

Since the duality between $\Delta$ and $\Delta^*$ establishes
a one-to-one correspondence between $1$-dimensional (respectively
$2$-dimensional) faces of $\Delta$ and $2$-dimensional (respectively
$1$-dimensional) faces of $\Delta^*$, we again obtain:

\begin{coro}
Let $\Delta$ be a reflexive polyhedron of dimension $4$, $\Delta^*$ the
dual reflexive polyhedron. Let $\hat{Z}_f$ be a MPCP-desingularization of
a $\Delta$-regular Calabi-Yau hypersurface $\overline{Z}_f$
in ${\bf P}_{\Delta}$, $\hat{Z}_g$ be a MPCP-desingularization of
a $\Delta^*$-regular Calabi-Yau hypersurface $\overline{Z}_g$
in ${\bf P}_{\Delta^*}$. Then the Euler characteristics of two Calabi-Yau
$3$-folds $\hat{Z}_f$ and $\hat{Z}_g$ satisfy the following relation
\[ e(\hat{Z}_f) = - e(\hat{Z}_g). \]
\label{coro.euler}
\end{coro}
\medskip

\section{Mirror symmetry}

\subsection{Mirror candidates for hypersurfaces of degree $n+1$ in
${\bf P}_n$}

\hspace*{\parindent}
Consider the  polyhedron $\Delta_n$ in $M_{\bf Q} \cong {\bf Q}^n$ defined
by inequalities
\[ x_1 + \ldots +  x_n \leq  1, \; x_i \geq -1 \; (1 \leq i \leq n). \]
Then $\Delta_n$ is a reflexive polyhedron, and
${\cal F}(\Delta)$ is the family of
all hypersurfaces of degree $n+1$ in ${\bf P}_n = {\bf P}_{\Delta_n}$.

The polyhedron $\Delta_n$ has $n+1$  $(n-1)$-dimensional faces whose
interiors contain exactly $n$ integral points. These $n(n+1)$
integral points form the root system of type $A_n$.
The number $l(\Delta_n)$ equals ${ 2n + 1 \choose n }$. Thus the dimension
of the moduli space of ${\cal F}({\Delta}_n)$ equals
\[ { 2n +1 \choose n } - (n+1)^2.\]

The dual polyhedron $\Delta^*_n$  has $n +1$ vertices
\[ u_1 =(1,0, \ldots, 0), \ldots , u_n = (0, \ldots, 0, 1),
u_{n+1} = (-1, \ldots , -1). \]
The corresponding toric Fano variety ${\bf P}_{\Delta_n^*}$ is a singular
toric hypersurface ${\bf H}_{n+1}$ of degree $n+1$ in ${\bf P}_{n+1}$ defined
by
the equation
\[  \prod_{i =1}^{n+1} u_i = u_0^{n+1},  \]
where $(u_0: \ldots : u_{n+1})$ are homogeneous coordinates in
${\bf P}_{n+1}$.

Since the simplex $\Delta_n$ is $(n+1)$-times multiple of $n$-dimensional
elementary
simplex of degree 1, the degree $d(\Delta_n)$ equals $(n+1)^n$. On the other
hand,
the dual simplex $\Delta_n^*$ is the union of $n+1$ elementary simplices of
degree 1, i.e., $d(\Delta_n^*) = n+1$.

There exists a finite morphism of degree $(n+1)^{n-1} =
d(\Delta_n)/ d(\Delta_n^*)$  of reflexive pairs
\[ \phi \; :\; (\Delta^*_n, N) \rightarrow (\Delta_n, M), \]
where  $\phi(u_{n+1}) = (-1, \ldots , -1) \in \Delta_n$ and

\[ \phi(u_i) = (-1, \ldots, \underbrace{n}_{i}, \ldots, -1 ) \in \Delta_n. \]

It is easy to see now that
\[  M / \phi(N) \cong ({\bf Z}/(n+1){\bf Z})^{n-1}. \]

Let $(v_0 : v_1:  \ldots : v_n)$ be the homogeneous coordinates on
${\bf P}_n$. The corresponding to $\phi$ \^etale  mapping of smooth
quasi-projective toric Fano  varieties
\[ \tilde {\phi} \;:\; {\bf P}_n^{[1]} \rightarrow
{\bf H}_{n+1}^{[1]} \]
has the following  representation in homogeneous coordinates
\[ (v_0 : v_1 : \dots : v_n ) \mapsto (\prod_{i= 0}^n v_i :
v_0^{n+1} : v_1^{n+1} : \dots : v_n^{n+1} ) =
(u_0 : u_1:  \ldots : u_{n+1} ).\]

A Calabi-Yau hypersurface $\overline{Z}_f$  in ${\bf H}_{n+1}$
has an  equation
\[ f(u)  = \sum_{i=0}^{n+1} a_i u_i = 0. \]

Using \ref{n.def}, it is easy to show that $h^{n-2,1}(Z_f) = 1$. We can also
describe the moduli space of the family ${\cal F}(\Delta_n^*)$.
Since  ${\bf H}_{n+1}$  is invariant under the action
of the $n$-dimensional torus
\[ {\bf T} = \{ {\bf t}=(t_1, \ldots, t_{n+1}) \in ({\bf C}^*)^{n+1} \mid
t_1  \cdots t_{n+1} = 1 \}, \]
the equation
\[ {\bf t}^*(f(u))  =  \sum_{i=0}^{n+1} a_i t_i u_i = 0, \]
defines an isomorphic to $\overline{Z}_f$ hypersurface
$\overline{Z}_{t^*(f)}$. Moreover, multiplying
all coefficients $\{a_i\}\; (0 \leq i \leq n+1)$
by the same non-zero complex number $t_0$, we get also a ${\bf C}^*$-action.
Thus, up to the action of the $(n+1)$-dimensional torus
${\bf C}^* \times {\bf T}$ on $n+2$ coefficients $\{ a_i\}$, we get
the one-parameter mirror family of Calabi-Yau hypersurfaces in
${\bf H}_{n+1}$ defined by the equation
\begin{equation}
 f_{a}(u)  = \sum_{i =1}^d u_i  + a u_0 = 0,
 \label{one.mir}
\end{equation}
where the number
\[ a  =  a_0 (\prod_{j =1}^{n+1} a_j)^{\frac{-1}{n+1}} \]
is uniquely defined up to an $(n +1)$-th root of unity.

Using  the homogeneous coordinates $\{ v_i \}$ on ${\bf P}_n$, we can transform
this equation to the form
\begin{equation}
 {\tilde {\phi}}^*f_{a} = \sum_{i =0}^n v_i^{n+1}   +
a \prod_{i=0}^n v_i,
\label{equ.mir}
\end{equation}
where $a^{n+1}$ is a canonical parameter of the corresponding subfamily
in ${\cal F}(\Delta_n)$ of smooth $(n-1)$-dimensional hypersurfaces
in ${\bf P}_n$.
\bigskip

\begin{theo}
The candidate for mirrors relative to the
family of all smooth
Calabi-Yau hypersurfaces of degree $n+1$ in ${\bf P}_n$ is the one-parameter
family ${\cal F}(\Delta^*_n)$ of Calabi-Yau varieties consisting  of quotients
by the action of
the finite abelian group $({\bf Z}/ (n+1){\bf Z})^{n-1}$  of the
hypersurfaces defined by the equation $(\ref{equ.mir})$.
\label{first.ex}
\end{theo}

\begin{exam}
{\rm If we take  $n =4$, then the corresponding finite abelian group
is isomorphic to  $({\bf Z}/5{\bf Z})^3$ and we come to the mirror
family for the family of all
$3$-dimensional quintics in ${\bf P}_4$ considered in \cite{cand2}.}
\end{exam}
\bigskip

\subsection{A category of reflexive pairs}

\hspace*{\parindent}

The set of all reflexive pairs of dimension $n$ forms a  category ${\cal C}_n$
whose  morphisms are finite morphisms of reflexive pairs (\ref{morph}).
The correspondence between dual reflexive pairs
defines an involutive functor
\[ {\rm Mir}\; : \; {\cal C}_n \rightarrow {\cal C}_n^* \]
\[ {\rm Mir}(\Delta, M) = (\Delta^*, N) \]
which is an isomorphism of the category ${\cal C}_n$ with the dual
category ${\cal C}_n^*$.
\bigskip

It is natural to describe in ${\cal C}_n$  some morphisms  satisfying
universal properties.

\begin{opr}
{\rm Let
\[ \phi_0 \; :\; (\Delta_0 , M_0) \rightarrow (\Delta, M) \]
be a finite morphism of reflexive pair. The morphism $\phi_0$ is said to
be {\em minimal} if for any finite morphism
\[ \psi \; :\; (\Delta' , M') \rightarrow (\Delta, M) \]
there exists the  unique morphism
\[ \phi \; :\; (\Delta_0 , M_0) \rightarrow (\Delta', M') \]
such that $\phi_0 = \psi \circ \phi$.

A reflexive pair in ${\cal C}_n$ is called {\em a minimal reflexive pair},
if the identity morphism of this pair  is minimal.}
\end{opr}
\bigskip

Note the following simple property.

\begin{prop}
Assume that there exist two finite morphisms
\[ \phi_1  \; :\; (\Delta_1, M_1) \rightarrow (\Delta_2, M_2) \]
and
\[ \phi_2  \; :\; (\Delta_2, M_2) \rightarrow (\Delta_1, M_1). \]
Then $\phi_1$ and $\phi_2$ are isomorphisms.
\label{isom}
\end{prop}

\proof The degrees  of $\phi_1$ and $\phi_2$ are positive
integers. On the other hand, by \ref{etale},
\[ d(\phi_1) = d_{M_2}(\Delta_2)/ d_{M_1}(\Delta_1), \]
\[ d(\phi_2) = d_{M_1}(\Delta_1)/ d_{M_2}(\Delta_2). \]
Therefore, $d_{M_1}(\Delta_1) = d_{M_2}(\Delta_2) =1$, i.e., $\phi_1$ and
$\phi_2$ are isomorphisms of reflexive pairs.
\bigskip

\begin{coro}
Assume that
\[ \phi_0 \; :\; (\Delta_0 , M_0) \rightarrow (\Delta, M) \]
is a minimal morphism. Then $(\Delta_0, M_0)$ is a minimal reflexive pair.
\end{coro}

\begin{coro}
For any reflexive pair $(\Delta, M)$ there exists up to an isomorphism
at most one minimal pair $(\Delta_0, M_0)$ with a minimal morphism
\[ \phi_0 \; : \; (\Delta_0 , M_0) \rightarrow (\Delta, M). \]
\end{coro}
\bigskip

The next proposition completely describes minimal reflexive pairs and the set
of all finite morphisms to a fixed reflexive pair $(\Delta, M)$.
\medskip

\begin{prop}
Let $(\Delta, M)$ be a reflexive pair. Denote by $M_{\Delta}$ the sublattice
in $M$ generated by vertices of $\Delta$. Let $M'$ be an integral lattice
satisfying the condition
$M_{\Delta} \subset M' \subset M$.  Then
$(\Delta, M')$ is also  a reflexive pair, and
\[ (\Delta, M_{\Delta}) \rightarrow (\Delta, M' ) \]
is a minimal morphism.
\label{min.mor}
\end{prop}

The proof immediately  follows from the definition of reflexive pair
\ref{inver.p}.
\medskip

\begin{coro}
All reflexive pairs having a finite morphism to a fixed reflexive
pair $(\Delta, M)$ are isomorphic to $(\Delta, M')$ for some  lattice
$M'$ such that $M_{\Delta} \subset M' \subset M$.
\label{min.pair}
\end{coro}
\bigskip

\begin{opr}
{\rm Let
\[ \phi^0 \; :\; (\Delta , M) \rightarrow (\Delta^0, M^0) \]
be a finite morphism of reflexive pair. The morphism $\phi^0$ is said to
be {\em maximal} if for any finite morphism
\[ \psi \; :\; (\Delta , M) \rightarrow (\Delta', M') \]
there exists the  unique morphism
\[ \phi \; :\; (\Delta' , M') \rightarrow (\Delta^0, M^0) \]
such that $\phi^0 = \phi \circ \psi$.

A reflexive pair in ${\cal C}_n$ is called {\em a maximal reflexive pair},
if the identity morphism of this pair  is maximal.}
\end{opr}
\bigskip

If we apply the functor ${\rm Mir}$, we get from \ref{min.mor} and
\ref{min.pair} the following properties
of maximal reflexive pairs.

\begin{prop}
Let $(\Delta, M)$ be a reflexive pair. Then  there exists up to an isomorphism
the unique maximal reflexive pair $(\Delta, M^{\Delta})$ having  a maximal
morphism
\[ \phi \; :\; (\Delta , M) \rightarrow (\Delta, M^{\Delta}). \]
Moreover, the pair  $(\Delta, M^{\Delta}) $ is dual to the minimal  pair
$(\Delta^*, N_{\Delta^*})$ having the  morphism
\[ \phi^*\; : \; (\Delta^* , N_{\Delta^*}) \rightarrow (\Delta^*, N) \]
as minimal.
\end{prop}
\bigskip

\begin{coro}
All reflexive pairs having finite morphisms  from a fixed reflexive
pair $(\Delta, M)$ are isomorphic to $(\Delta, M')$ for some  lattice
$M'$ such that $M \subset M' \subset M^{\Delta}$.
\end{coro}
\bigskip

\begin{exam}
{\rm Let $(\Delta_n, M)$ and $(\Delta^*_n, N)$ be two reflexive pairs
from the previous section. Since the lattice $N$ is generated by
vertices of $\Delta_n^*$, the reflexive pair  $(\Delta_n^*, N)$ is minimal.
Therefore, $(\Delta_n, M)$ is a maximal reflexive pair. }
\end{exam}
\bigskip

\subsection{A Galois correspondence}

\hspace*{\parindent}

The existence of a finite morphism of reflexive pairs
\[ \phi\; : \; (\Delta_1, M_1) \rightarrow (\Delta_2, M_2) \]
implies the following main geometric relation between
Calabi-Yau hypersurfaces in toric
Fano varieties ${\bf P}_{\Delta_1, M_1}$ and ${\bf P}_{\Delta_2, M_2}$.
\medskip

\begin{theo}
The Calabi-Yau hypersurfaces in ${\bf P}_{\Delta_1, M_1}$ are quotients of
some Calabi-Yau hypersurfaces in ${\bf P}_{\Delta_2, M_2}$ by the action of
the dual to $M_2/ \phi(M_1)$ finite abelian group.
\label{quot}
\end{theo}

\proof Consider the dual finite morphism
\[ \phi^* \; : \; (\Delta_2^*, N_2) \rightarrow (\Delta_1^*, N_1). \]
By \ref{etale}, $\phi^*$ induces a finite  \^etale morphism
\[ \tilde {\phi}^* \; : \;  {\bf P}_{\Delta_2,M_2}^{[1]}
\rightarrow  {\bf P}_{\Delta_1,M_1}^{[1]}  \]
of smooth quasi-projective toric Fano varieties.  This morphism is defined
by the surjective homomorphism of $n$-dimensional algebraic tori
\[ {\gamma}_{\phi} \; :\; {\bf T}_{\Delta_2}
\rightarrow {\bf T}_{\Delta_1} \]
whose kernel is dual to the cokernel of the homomorphism of the groups of
characters
\[ \phi \; :\; M_1 \rightarrow  M_2.\]
Therefore, $\tilde {\phi}^*$ is the quotient by the action of the
finite abelian group  $(M_2/ \phi(M_1))^* =  N_1 / \phi^*(N_2)$.
The pullback of the anticanonical
class of  ${\bf P}_{\Delta_1,M_1}^{[1]}$
is the anticanonical class of
${\bf P}_{\Delta_2,M_2}^{[1]}$.
Applying \ref{galois}, we obtain that
the smooth quasi-projective Calabi-Yau hypersurfaces
${\hat Z}_{f,\Delta_1}$ are \^etale quotients by $N_1 / \phi^*(N_2)$
of some smooth
quasi-projective Calabi-Yau hypersurfaces
${\hat Z}_{\tilde {\phi}^*f,\Delta_2}$.
\bigskip

\begin{coro}
The mirror mapping for families of Calabi-Yau hypersurfaces in toric
varieties satisfies the following Galois correspondence:

If a family ${\cal F}(\Delta_1)$ is a quotient of a family
${\cal F}(\Delta_2 )$ by a finite abelian group ${\cal A}$, then
the mirror family ${\cal F}(\Delta_2^*)$ is a quotient of the mirror
family ${\cal F}(\Delta_1^*)$ by the dual finite abelian
group ${\cal A}^*$.
\end{coro}
\bigskip

\begin{opr}
{\rm
Let $(\Delta, M)$ be a reflexive pair,
$(\Delta^*,  N)$  the dual reflexive pair.
Denote  by  $N_{\Delta^*}$
the sublattice  in $N$  generated by  vertices of  $\Delta^*$.
The finite abelian  group $\pi_1 (\Delta, M) = N/ N_{\Delta^*}$
 is called  the {\em fundamental group of the pair } $(\Delta, M)$. }
\end{opr}
\bigskip

The fundamental group $\pi_1 (\Delta, M)$ defines  a contravariant
functor from the category ${\cal C}_n$ to the category of finite
abelian groups with injective homomorphisms.
\bigskip

\begin{prop}
Let $(\Delta, M)$ be a reflexive pair, ${\bf P}_{\Delta, M}$
the corresponding toric Fano variety. Then the fundamental group
$\pi_1 (\Delta, M)$ is
isomorphic to the algebraic $($and topological$)$ fundamental group
\[ \pi_1 ({\bf P}^{[1]}_{\Delta} ). \]
In particular, the reflexive pair $(\Delta, M)$ is maximal if and only if
${\bf P}^{[1]}_{\Delta}$ is simply connected.
\end{prop}

\proof The statement immediately
follows from the description of the fundamental
group of toric varieties in \ref{fund.group}.
\bigskip

\begin{opr}
{\rm Let $(\Delta, M)$ be a reflexive pair, $(\Delta^*, N)$ the dual
reflexive pair, $(\Delta, M_{\Delta})$ and $(\Delta^*, N_{\Delta^*})$ are
minimal pairs, $(\Delta, M^{\Delta})$ and $(\Delta, N^{\Delta^*})$
are maximal pairs. The quotients
\[ \pi_1 (\Delta) = N^{\Delta^*} / N_{\Delta^*} \; {\rm and }\;
\pi_1 (\Delta^*) = M^{\Delta} / M_{\Delta } \]
 is called the {\em fundamental groups} of the reflexive polyhedra $\Delta$
and  $\Delta^*$ respectively. }
\end{opr}
\medskip

It is clear, $\pi_1 (\Delta)$ and $\pi_1 (\Delta^*)$ are isomorphic
dual finite abelian  groups.
\medskip

\begin{opr}
{\rm Assume that for a reflexive pair $(\Delta, M)$
there exists an isomorphism between two maximal
reflexive pairs
\[ \phi\; :\; (\Delta, M^{\Delta}) \rightarrow (\Delta^*, N^{\Delta^*}) .\]
Then we call  $\Delta$ {\em a selfdual reflexive polyhedron}. }
\label{selfdual}
\end{opr}
\medskip

If $\Delta$ is selfdual, then $\Delta$ and $\Delta^*$ must have the same
combinatorial type (see \ref{dual.edge}). By \ref{quot}, we obtain.

\begin{prop}
Let $(\Delta, M)$ be a reflexive pair such that $\Delta$ is selfdual.
Then ${\cal F}(\Delta)$ and ${\cal F}(\Delta^*)$ are
quotients respectively by $\pi_1(\Delta, M)$ and
$\pi_1(\Delta^*, N)$  of some  subfamilies in the family
of Calabi-Yau hypersurfaces
corresponding to two isomorphic maximal reflexive pairs
$(\Delta, M^{\Delta})$ and $(\Delta, N^{\Delta^*})$.
Moreover, the order of $\pi_1(\Delta)$ equals to the product of
oders of $\pi_1(\Delta, M)$ and $\pi_1(\Delta^*, N)$.
\end{prop}
\bigskip

\subsection{Reflexive simplices}

\hspace*{\parindent}

In this section we consider Calabi-Yau families ${\cal F}(\Delta)$,
where $\Delta$ is a reflexive simplex of dimension $n$.
Let $\{ p_0, \ldots, p_n \}$ be vertices of $\Delta$.
There exists the  unique linear relation among $\{ p_i \}$
\[ \sum_{i = 0}^n w_i p_i = 0, \]
where $w_i > 0$ ($0 \leq i \leq n$) are integers
and $g.c.d. (w_i ) = 1.$

\begin{opr}
{\rm The coefficients $w = \{ w_0, \cdots , w_n \}$ in the above linear
relations are called  the {\em weights of the reflexive simplex} $\Delta$.}
\end{opr}
\bigskip

Let $\Delta^*$ be the dual reflexive simplex,  ${l}_0 , \ldots , {l}_n$
vertices of $\Delta^*$, $p_0, \dots , p_n$ vertices of $\Delta$.
By definition
\ref{inver.p}, we may  assume that
\[ \langle p_i , {l}_j \rangle = -1 \; (i \neq j ), \]
i.e., that the equation $\langle x , {l}_j \rangle = -1$ defines the affine
hyperplane in ${M}_{\bf Q}$ generated by  the $(n-1)$-dimensional face of
$\Delta$ which does not contain  $p_j \in \Delta$.

\begin{opr}
{\rm The  $(n+1)\times(n+1)$-matrix with integral coefficients
\[ B(\Delta) = ( b_{ij} ) = ( \langle p_i, {l}_j \rangle ) \]
is called  {\em the matrix of the reflexive simplex $\Delta$}. }
\end{opr}
\bigskip

\begin{theo}
Let  $\Delta$ be a reflexive $n$-dimensional simplex. Then

{\rm (i)} the matrix $B(\Delta)$ is symmetric and its rank equals $n$;

{\rm (ii)} the diagonal coefficients $b_{ii}$ $( 0 \leq i \leq n)$
are positive and satisfy the equation
\[\sum_{i =0}^n \frac{1}{b_{ii} + 1} = 1;\]

{\rm (iii)} the weights $\{ w_0, \ldots , w_n \}$ of $\Delta$ are
the primitive integral solution of the linear  homogeneous
system with the matrix $B(\Delta)$. Moreover,
\[ w_i = \frac{l.c.m. (b_{ii} +1)}{b_{ii} +1}. \]
\label{ref.mat}
\end{theo}

\proof The statement (i) follows from the fact that
${\rm rk}\, \{ p_i \} = {\rm rk}\, \{ l_j \} = n $. One gets
(ii) by the direct computation of the determinant of $B(\Delta)$ as a function
on coefficients $b_{ii}$ ($0 \leq i \leq n$).  Finally,
(iii) follows from (ii) by checking  that
\[ \{ 1/(b_{11} +1), \dots , 1/(b_{nn} +1) \} \]
is a solution of the linear homogeneous system with the matrix $B(\Delta)$.
\bigskip

\begin{coro}
The matrix $B(\Delta)$ depends only on the weights of $\Delta$.
\end{coro}
\bigskip

Let $\Delta$ be a reflexive simplex with weights $w = \{ w_i \}$. Using the
vertices $\{ l_j \}$ of the dual reflexive simplex $\Delta^*$, we can
define the homomorphism
\[ \iota_{\Delta} \;:\;  M \rightarrow {\bf Z}^{n+1},\]
where
\[ \iota_{\Delta} (m) =
(\langle m, {l}_0 \rangle, \ldots , \langle m, {l}_n \rangle ). \]
Obviously, $\iota_{\Delta}$ is injective and the image of $\iota_{\Delta}$
is contained in the $n$-dimensional sublattice $M(w)$
in ${\bf Z}^{n+1}$ defined by the equation
\[ \sum_{i =0}^n  w_i x_i = 0.\]
Note that the image $\iota_{\Delta}(p_i)$ is the $i$-th row of $B(\Delta)$. We
denote by $\Delta(w)$ the convex hull  of the points
$\{ \iota_{\Delta}(p_i)\}$ in $M_{\bf Q}(w)$.
\medskip

\begin{theo}
The pair $(\Delta(w), M(w))$ is reflexive and satisfies the following
conditions$:$

{\rm (i)} the corresponding to $(\Delta(w), M(w))$
toric Fano variety ${\bf P}_{\Delta(w)}$
is the weighted projective space ${\bf P}(w_0, \cdots, w_n)$;

{\rm (ii)}
\[ \iota_{\Delta} \; : \; (\Delta, M) \rightarrow
(\Delta(w), M(w)) \]
is a finite  morphism of reflexive pairs;

{\rm (iii)} $(\Delta(w), M(w))$ is a  maximal reflexive pair.
\label{matrix}
\end{theo}

\proof The reflexivity of $(\Delta(w), M(w))$ and the condition (ii)
follow immediately from
the  definition of $\iota_{\Delta}$, since
\[ \Delta(w) = \{ (x_0, \ldots, x_n) \in {\bf Q}^{n+1} \mid
\sum_{i =0}^n w_i x_i =0, \; x_i \geq  -1 \; (0 \leq i \leq n) \}. \]

(i) The shifted by $(1, \ldots, 1)$ convex polyhedron
\[ \Delta^{(1)}(w) = \Delta(w) + (1, \ldots, 1) \]
is the intersection of ${\bf Q}^{n+1}_{\geq 0}$ and the affine hyperplane
\[ w_0 x_0 + \cdots + w_n x_n = w_0 + \cdots + w_n =
l.c.m.\{ b_{ii} +1 \} = d. \]
Therefore, the integral points in $\Delta^{(1)}(w)$ can be identified
with all possible monomials of degree $d$ in $n+1$ $w$-weighted
independent variables.

(iii) Assume that there exists a finite morphism
\[ \phi\; :\; (\Delta(w), M(w)) \rightarrow (\Delta', M')\]
of reflexive pairs. Obviously, $\Delta'$ must be also a  simplex.
Since  the linear mapping $\phi$    does  not change linear relations,
$\Delta$ and $\Delta'$  must have the same weights  $w = \{w_i \}$.
Therefore, by (ii), there exists a finite morphism
\[ \iota_{\Delta'} \; : \; (\Delta', M') \rightarrow
(\Delta(w), M(w)). \]
Therefore, by \ref{isom}, the reflexive pairs $(\Delta(w), M(w))$ and
$(\Delta', M')$ are isomorphic.
\bigskip

Since $B(\Delta) = B(\Delta^*)$ (see \ref{ref.mat} ), we obtain:

\begin{coro}
Any reflexive simplex  $\Delta$ is selfdual.
\label{dual.simp}
\end{coro}
\bigskip

\subsection{Quotients of Calabi-Yau hypersurfaces
in weighted projective spaces}

\hspace*{\parindent}

Theorem \ref{matrix} implies

\begin{coro}
The fundamental group $\pi_1(\Delta)$ of a  reflexive simplex
$\Delta$ depends only on the  weights $w = \{ w_i \}$. This group is
dual to the quotient of the lattice $M(w) \subset {\bf Z}^{n+1}$ by the
sublattice $M_B(\Delta)$ generated by rows of the matrix $B(\Delta)$.
\end{coro}

Now we calculate  the group  $M(w) / M_B(w)$ explicitly.
\bigskip

Consider two integral sublattices of
rang $n+1$ in ${\bf Z}^{n+1}$
\[ \tilde{M} (w) = M(w) \oplus {\bf Z}\langle (1, \ldots, 1) \rangle, \]
and
\[ \tilde{M}_B (w) = M_B(w) \oplus {\bf Z}\langle (1, \ldots, 1) \rangle. \]
Note that $\tilde{M}(w) / \tilde{M}_B(w) \cong M(w) / M_B(w)$.

Let $\mu_r$ denote the group of complex $r$-th roots of unity. Put
$d_i = b_{ii} +1$ ($ 0 \leq i \leq n$), $d = l.c.m.\{ d_i \}$.

The sublattice $\tilde{M}(w)$ is the kernel of the surjective
homomorphism
\[ \gamma_w \;: \; {\bf Z}^{n+1} \rightarrow \mu_d, \]
\[ \gamma_w (a_0, \ldots , a_n) = g^{w_0a_0 + \cdots + w_n a_n}, \]
where $g$ is a generator of $\mu_d$.

The sublattice $\tilde{M}_B(w)$ is generated by $(1, \ldots, 1)$ and
\[ (d_0, 0, \ldots, 0), (0, d_1, \ldots, 0),
\dots, (0, \ldots, 0, d_n). \]

Therefore, $\tilde{M}_B(w)$ can be represented as the sum of the infinite
cyclic group  generated by $(1, \ldots, 1)$ and
the kernel of the surjective homomorphism
\[ \gamma \;: \; {\bf Z}^{n+1} \rightarrow \mu_{d_0} \times \mu_{d_1}
\times \cdots \times \mu_{d_n} ,  \]
\[ \gamma (a_0, \ldots , a_n) = g_0^{a_0} g_1^{a_1} \cdots g_n^{a_n}, \]
where $g_i$ is a generator of $\mu_{d_i}$ ($ 0 \leq i \leq n$).
The order of the element $(1, \ldots, 1)$ modulo ${\rm ker}\, \gamma$
equals  $d$.
Thus,
\[ {\bf Z}^{n+1} /  \tilde{M}_B(w) \cong (\mu_{d_0} \times \mu_{d_1}
\times \cdots \times \mu_{d_n}) / \mu_d ,\]
where the subgroup $\mu_d$ is generated by
$g_0g_1 \cdots g_n$. Finally, we obtain.

\begin{theo}
The fundamental group $\pi_1(\Delta^*) = M(w)/M_B(w)$ of the
reflexive simplex $\Delta^*$ with weights
$w = \{ w_i \}$ is isomorphic  to the  kernel
of the surjective homomorphism
\[ \overline{\gamma}_w  \; : \; (\mu_{d_0} \times \mu_{d_1}
\times \cdots \times \mu_{d_n}) / \mu_d  \rightarrow \mu_d, \]
\[ \overline{\gamma}_w ( g_0^{a_0} g_1^{a_1} \cdots g_n^{a_n} ) =
g^{w_0a_0 + \cdots + w_n a_n}. \]
\label{fund1}
\end{theo}

By duality, we conclude.

\begin{coro}
The fundamental group $\pi_1(\Delta)$ of a reflexive simplex $\Delta$ with
weights $w$ is  isomorphic to the kernel of the surjective
homomorphism
\[ (\mu_{d_0} \times \mu_{d_1} \times \cdots \times \mu_{d_n}) / \mu_d
\rightarrow \mu_d, \]
where the homomorphism to $\mu_d$ is the product of complex numbers
in $\mu_{d_0}, \mu_{d_1}, \ldots , \mu_{d_n}$, and the
embedding of $\mu_d$ in $\mu_{d_0} \times \mu_{d_1}
\times \cdots \times \mu_{d_n}$ is defined by
\[ g \mapsto ( g^{w_0}, \ldots, g^{w_n}). \]
\label{fund2}
\end{coro}

\begin{coro}
The order  of $\pi_1(\Delta)$ in the above theorem equals
\[ \frac{d_0 d_1 \cdots d_n}{d^2}. \]
\end{coro}
\bigskip

\begin{exam}
{\rm Let $(d_0, d_1, \ldots , d_n )$ $(d_i > 0)$ be an
integral solution of the equation
\begin{equation}
\label{weights.proj}
 \sum_{i =0}^{n} \frac{1}{d_i} = 1.
\end{equation}
Then the quasi-homogeneous equation
\[ v_0^{d_0} + v_1^{d_1} + \cdots + v_n^{d_n} = 0\]
defines a $\Delta(w)$-regular Calabi-Yau hypersurface of Fermat-type in the
weighted projective space
\[ {\bf P}_{\Delta(w)} = {\bf P}(w_0, \ldots, w_n) , \]
where
\begin{equation}
 w_i = \frac{l.c.m. (d_i)}{d_i}.
\label{new.weights}
\end{equation}}
\end{exam}

\begin{coro}
The family ${\cal F}(\Delta(w))$ of
Calabi-Yau hypersurfaces in the weighted projective
space ${\bf P}_{\Delta(w)}$ consists of deformations of Fermat-type
hypersurfaces.
If $\Delta$ is a reflexive  simplex  with weights $w =\{ w_i \}$, then
the corresponding family ${\cal F}(\Delta)$ in ${\bf P}_{\Delta}$
consists of quotients of  some subfamily in ${\cal F}(\Delta(w))$
by the  action of the finite abelian group $\pi_1( \Delta,M)$.
\end{coro}

If we consider a special case $n =4$, we obtain as a corollary
the result of Roan in \cite{roan1}. To prove this, one should
use our general result on Calabi-Yau $3$-folds
constructed from $4$-dimensional reflexive polyhedra (see \ref{euler.iso})
and apply the following simple statement.

\begin{prop} Let $(\Delta, M)$ be a reflexive pair such that
$\Delta$ is a $4$-dimensional reflexive simplex with weights
$(w_0, \ldots, w_4)$. Then
the family ${\cal F}(\Delta)$ consists of quotients by $\pi_1(\Delta, M)$
of Calabi-Yau hypersurfaces in the weighted projective space
${\bf P}(w_0, \ldots, w_4)$ whose equations are invariant under the
canonical diagonal action of $\pi_1(\Delta, M)$ on
${\bf P}(w_0, \ldots, w_4)$.
\end{prop}

\end{document}